\magnification=\magstep1
\vbadness=10000
\hbadness=10000
\tolerance=10000

\def\Aut{{\rm Aut}}
\def\height{{\rm height}}

\proclaim Classification of positive definite lattices. 

Richard E. Borcherds, 
\footnote{$^*$}{ Supported by 
NSF grant DMS-9970611.}

Mathematics department,
Evans Hall \#3840,
University of California at Berkeley,
CA 94720-3840
U.S.A.

e-mail: reb@math.berkeley.edu

www home page  www.math.berkeley.edu/\hbox{\~{}}reb

\bigskip

\proclaim Contents.

1. An algorithm for classifying vectors in some Lorentzian lattices. 

2. Vectors in the lattice $II_{1,25}$.

3. Lattices with no roots. 

Table 0:  Primitive norm $0$ vectors in $II_{1,25}$.

Table 1: Norm $2$ vectors in $II_{1,25}$.

Table 2: Norm $4$ vectors in $II_{1,25}$.

\proclaim 1.~Classification of positive  norm vectors. 

In this paper we describe an algorithm for classifying orbits of
vectors in Lorentzian lattices. The main point of this is that
isomorphism classes of positive definite lattices in some genus often
correspond to orbits of vectors in some Lorentzian lattice, so we can
classify some positive definite lattices. 
Section 1 gives an overview of this algorithm, and in section 2 we describe
this algorithm more precisely for the case of $II_{1,25}$, and as an
application we give the classification of the 665 25-dimensional
unimodular positive definite lattices and the 121 even 25 dimensional
positive definite lattices of determinant 2 (see tables 1 and 2). In section 3
we use this algorithm to  show that there is a unique 26
dimensional unimodular positive definite lattice with no roots.
Most of the results of this paper are taken from the unpublished
manuscript [B], which contains more details and examples.
For general facts about lattices used in this paper see [C-S], especially
chapters 15--18 and  23--28. 

Some previous enumerations of unimodular lattices include Kneser's
list of the unimodular lattices of dimension at most 16 [K], Conway
and Sloane's extension of this to dimensions at most 23 [C-S chapter
16], and Niemeier's enumeration [N] of the even 24 dimensional
ones. All of these used some variation of Kneser's neighborhood
method [K], but this becomes very hard to use for odd lattices of
dimension 24, and seems impractical for dimension at least 25 (at
least for hand calculations; computers could probably push this
further). The method used in this paper works well up to 25
dimensions, could be pushed to work for 26 dimensions, and does not
seem to work at all beyond this.

We use the ``$(+,-,-,\cdots,-)$'' sign convention for Lorentzian
lattices $L$, so that the reflection we are interested in are
(usually) those of {\it negative} norm vectors of $L$. We fix one of
the two cones of positive norm vectors and call it the positive
cone. The norm 1 vectors in the positive cone form a copy of
hyperbolic space in the usual way.  We assume that we are given a
group $G$ of automorphisms of a Lorentzian lattice $L$, such that $G$
is the semidirect product of a normal subgroup $R$ generated by
reflections of some negative norm vectors, and a group $\Aut(D)$ of
automorphisms preserving a fundamental domain domain $D$ of $R$ in
hyperbolic space.  We assume that all elements of $u\in L$ having
non-negative inner product with all simple roots of $R$ have norm
$(u,u)$ at least 0 (this is just to eliminate some degenerate cases).
If $L$ is a lattice then $L(-1)$ is the lattice $L$ with all norms
multiplied by $-1$.
We  use Conway's convention of using small letters $a_n$, $d_n$,
$e_n$ for the spherical Dynkin diagrams, and capital letters $A_n$,
$D_n$, $E_n$ for the corresponding affine Dynkin diagrams.
The Weyl vector of a root system is the vector $\rho$ such that 
$(\rho,r)=-r^2/2$ for any simple root $r$. 

We want to find the orbits of positive norm vectors of the positive
cone of $L$ under the group $G$.  Every positive norm vector of the
positive cone of $L$ is conjugate under $R$ to a unique vector in $D$,
so it is enough to classify orbits of vectors $u$ in $D$ under
$\Aut(D)$.

The algorithm works by trying to reduce a vector $u$ of $D$ to a
vector of smaller norm by adding a root of $u^\perp$ to $u$. There are
three possible cases we need to consider:
\item{(1)} 
There are no roots in $u^\perp$. 
\item{(2)} 
There is a root $r$ in $u^\perp$ such that $u+r\in D$. 
\item{(3)} 
There is at least one root in $u^\perp$, but if $r$ is a root in
$u^\perp$ then $u+r$ is never in $D$.

We try to deal with these three cases as follows. 

If there are no roots in $u^\perp$, then we assume that $D$ contains a
non-zero vector $w$ such that $(r,w)\le (r,u)$ for any simple root $r$
and any vector $u\in L$ in the interior of $D$.  Then $u-w$ has inner
product at least 0 with all simple roots, so it also lies in $D$ and
has smaller norm than $u$ unless $u$ is a multiple of $w$ and
$w^2=0$. So we can reduce $u$ to a vector of smaller norm in $D$.  The
existence of a vector $w$ with these properties is a very strong
condition on the lattice $L$.

\noindent 
{\bf Example 1.1.} The lattices $II_{1,9}$ and $II_{1,17}$ have
properties 1 and 2; this follows easily from Vinberg's description [V85] of
their automorphism groups. Conway showed that the lattice $II_{1,25}$
also has these properties; see the next section. The lattices
$II_{1,8n+1}$ for $n\ge 4$ do not have these properties; but the
Minkowski-Siegel mass formula shows that these lattices have such vast
numbers of orbits of positive norm vectors that there seems little
point in classifying them.

\noindent
{\bf Example 1.2.} It follows from [B90] that several lattices that
are fixed points of finite groups acting on $II_{1,25}$ also have
a suitable vector $w$.  For example the lattice $II_{1,1}\oplus BW(-1)$, where
$BW$ is the Barnes-Wall lattice, has this property. Some of the norm 0
vectors correspond to the 24 lattices in the genus of $BW$ classified
in [S-V]; the remaining orbits of norm 0 vectors should not be hard to
find.

\noindent
{\bf Example 1.3.}  Take $L$ to be the lattice $I_{1,9}$ and $R$ to be
the group generated by reflections of norm $-1$ vectors. (This has
infinite index in the full reflection group.) Then the lattice has a
Weyl vector for the reflection group as in [B90], so we can apply the
algorithm to this reflection group. (However it is not entirely clear
what the point of doing this is, as it is easier to use the full
reflection group of the lattice!)

Next we look at the second case when $u^\perp$ has a root $r$ such
that $v=u-r$ is in $D$. Then $-r$ is in the fundamental domain of the
finite reflection group of $u^\perp$, so $r$ is a sum of the simple
roots of $u^\perp$ with the usual multiplicities.

For $u$ in $D$ we let $R_i(u)$ be the simple roots $u$ of $D$ that
have inner product $i(r,r)/2$ with $u$, so $R_i(u)$ is empty for $i<0$
and $R_0(u)$ is the Dynkin diagram of $u^\perp$. We write $S(u)$ for
$R_0(u)\cup R_1(u)\cup R_2(u)$. Then given $S(v)$ we can find all
vectors $u$ of $D$ that come from $v$ as in (3) above, and $S(u)$ is
contained in $S(v)$. By keeping track of the action of $\Aut(D,v)$ on
$S(v)$ for vectors $v$ of $D$ we can find all possible vectors $v$
constructed in this way from $v$, together with the sets $S(u)$.

Finally, the third case, when there is at least one root in $u^\perp$,
but if $r$ is a root in $u^\perp$ then $v-r$ is never in $D$, has to
be dealt with separately for each lattice $L$. In practice it does not
present too much difficulty for lattices with a vector $w$ as in case
1. See the next section for the example of $L=II_{1,25}$.

The following two lemmas will be used later to prove some properties
of the root systems of 25 dimensional lattices.

\proclaim Lemma 1.4.
Suppose that reflection in $u^\perp$ is an automorphism of $L$.  Then
there is an automorphism $\sigma$ of $L$ (of order 1 or 2) with the
following properties:
\item{(1)} $\sigma$ fixes $D$.
\item{(2)} If $\sigma$ fixes $w$, then $w$ is a linear combination of
$u$ and the roots of $L$ in $u^\perp$.

Proof. There is an automorphism of $L$ acting as 1 on $u$ and as $-1$
on $u^\perp$, given by the product of $-1$ and reflection in
$u^\perp$.  As this automorphism fixes $u\in D$, we can multiply it by
some (unique) element of the reflection group of $u^\perp$ so that the
product $\sigma$ fixes $D$.  The element $\sigma$ acts as $-1$ on the
space orthogonal to $u$, $z$, and all roots of $R$ in $u^\perp$, which
implies assertion (2) of the lemma 1.4.

\proclaim Lemma 1.5. 
Suppose that there is a norm 0 vector $z$ such that $(z,u)=2$, where
$u$ is a vector in $D$. Then there is an automorphism $\sigma$ of $L$
with the following properties:
\item{(1)} 
$\sigma$ fixes $D$.
\item{(2)} 
If $\sigma$ fixes $w$, then $w$ is a linear combination of $u$, $z$,
and the roots of $L$ in $u^\perp$.

Proof. If $M$ is the lattice spanned by $z$ and $u$ then $M$ has the
property that all elements of $M'/M$ have order 1 or 2. So there is an
automorphism of $L$ acting as 1 on $M$ and $-1$ on $M^\perp$.  The
result now follows as in the proof of lemma 1.4.  This proves lemma
1.5.

{\bf Remark.} It is usually easy to classify all orbits of negative
norm vectors $u$ in Lorentzian lattices, because this is closely
related to the classification of the {\it indefinite} lattices
$u^\perp$, and by Eichler's theorem [E] indefinite lattices in dimension
at least 3 are classified by the spinor genus (which in practice is
often determined by the genus). For example, it is easy to give a
proof along these lines that if $n>0$ and $m>0$ then $II_{1,8n+1}$ has
a unique orbit of primitive vectors of norm $-2m$.

\proclaim 2.~Vectors in the lattice $II_{1,25}$.

In this section we specialize the algorithm of the previous section to
the lattice $II_{1,25}$.

Note that orbits norm 4 vectors $u$ of $II_{1,25}$ correspond
naturally to 25 dimensional positive definite unimodular lattices,
because $u^\perp$ is isomorphic to the lattice of even vectors in a 25
dimensional unimodular negative definite lattice. In particular we can
classify the 665 positive definite 25 dimensional unimodular lattices,
as in table 2; this is the main application of the algorithm of the
previous section. Similarly norm 2 vectors of $II_{1,25}$ correspond
to 25 dimensional even positive definite lattices of determinant 2.

First we have to show the existence of a vector $w$ satisfying the
property of section 1. This follows from Conway's theorem [C85]
stating that the reflection group of $II_{1,25}$ has a Weyl vector $w$
of norm 0, with the property that $(w,r)=1$ for all simple roots $r$
of the reflection group.  Conway's proof depends on the rather hard
classification of the ``deep holes'' in the Leech lattice in [C-P-S];
there is a proof avoiding these long calculations in [B85].  It seems
likely that 26 is the largest possible dimension of a lattice with a
suitable vector $w$.

Next we have to classify the vectors $u$ of $D$ such that $u^\perp$
has roots but $u+r$ is not in $D$ for any root $r\in u^\perp$.  One
obvious way this can happen is if $u$ has norm 0, so we have to
classify the norm 0 vectors in $II_{1,25}$. In any lattice
$L=II_{8n+1,1}$ the orbits of primitive norm 0 vectors $z$ correspond
to the $8n$-dimensional even negative definite unimodular lattices
$z^\perp/z$. So the orbits of primitive norm 0 vectors of $II_{1,25}$
correspond to the 24 Niemeier lattices ([C-S]). The non-primitive norm
0 vectors are of course either 0 or a positive integer multiple of a
primitive norm 0 vector, so this gives the classification of all
orbits of norm 0 vectors in $II_{1,25}$; see table 0.

Next suppose that $u$ is a positive norm vector of $D$ with $(u,u)=2n$
and $r$ is a highest root in $u^\perp$ such that $urr$ is not in $D$.
Then $u-r$ is conjugate under the reflection group to some vector $v$
such that $(v,u)<(u-r,u)$. But $(v,u)^2\ge (u,u)(v,v) =2n(2n-2)$ and
$(v,u)<(u-r,u)=2n$, so $(v,u)=2n-1$.  So if $z=u-v$ then $(z,u)=1$ and
$z^2=0$.  If we put $z'=u-nz$ then $z$ and $z'$ are norm 0 vectors
with $(z,z')=1$ and $u=nz+z'$. So $II_{1,25}=B\oplus \langle
z,z'\rangle$ for some Niemeier lattice $B$. If this Niemeier lattice
has roots, then adding some of these roots to $r$ gives a vector in
$D$ by the previous argument, so $B$ must be the Leech lattice so we
can assume that $z$ is in the orbit of $w$. If $n>1$ then there are no
roots in $u^\perp$, and if $n<1$ then $(u,u)\le 0$, so we must have
$n=1$.  So the only possibility for $u$ is that it is a norm 2 vector
in the orbit of $w+w'=2w+r$, where $r$ is a simple root.

Putting everything together gives the following list of the vectors
$u\in D$ such that $u^\perp$ has roots but $u-r$ is not in $D$ for any
root $r\in u^\perp$:
\item{1.} 
The zero vector.
\item{2.} 
The norm 0 vectors $nz$ for $n\ge 1$ and $z$ a primitive norm 0 vector of $D$
corresponding to some Niemeier lattice other than the Leech lattice.
The vectors for a given Niemeier lattice and a given value of $n$ are
all conjugate under $\Aut(D)$.
\item{3.}
The norm 2 vectors of the form $2w+r$ for a simple root $r$ of
$D$. These form one orbit under $\Aut(D)$.

\proclaim Lemma 2.1.
Suppose $u,v\in D$, $u^2=2n$, $v^2=2(n-1)$ and $(v,u)=2n$. Then 
$$\eqalign{
R_0(u)&\subseteq R_0(v)\cup R_1(v)\cup R_2(v)=S(v)\cr
R_i(u)&\subseteq R_0(v)\cup R_1(v)\cup \cdots\cup R_i(v)\hbox{ for } i\ge 1.\cr
}$$

Proof. The vector $v$ is in $D$, so $v=u+r$ for some highest root $r$
of $u^\perp$. The vector $r$ has inner product 0, $1$, or $2$ with all
simple roots of $u^\perp$, and $-r$ is a sum of roots of $R_0(u)$ with
positive coefficients, so $r$ has inner product $\ge 0$ with all
simple roots of $D$ not in $R_0(u)$. The lemma follows from this and
the fact that $(v,s)=(u,s)+(r,s)$ for any simple root $s$ of $D$. This
proves lemma 2.1.

We now start with a vector $v$ of norm $2(n-1)$ and try to reconstruct
$u$ from it. The vector $u-v$ is a highest root of some component of
$R_0(u)$, and $R_0(u)$ is contained in $S(v)$, so we should be able to
find $u$ from $S(v)$. By lemma 2.1 $S(u)$ is contained in $S(v)$, so
we can repeat this process with $u$ instead of $v$.  The following
theorem shows how to construct all possible vectors $u$ as in lemma
2.1 from $v$ and $S(v)$.

\proclaim Theorem 2.2.
Suppose that $v$ has norm $2(n-1)$ and is in $D$ (so $n\ge 1$). Then
there are bijections between
\item{(1)}
Norm $2n$ vectors $u$ of $D$ with $(u,v)=2n$. 
\item{(2)}
Simple spherical Dynkin diagrams $C$ contained in the Dynkin diagram
$\Lambda$ of $D$ such that if $r$ is the highest root of $C$ and $c$
in $C$ satisfies $(c,r)=i$, then $c$ is in $R_i(v)$.
\item{(3)}
Dynkin diagrams $C$ satisfying one of the following three conditions: 
\itemitem{Either}
$C$ is an $a_1$ and is contained in $R_2(v)$, 
\itemitem{or}
$C$ is an $a_n$ ($n\ge 2$) and the two endpoints of $C$ are in
$R_1(v)$ while the other points of $C$ are in $R_0(v)$,
\itemitem{or}
$C$ is $d_n$ ($n\ge 4$), $e_6$, $e_7$, or $e_8$ and the unique point of
$C$ that has inner product $1$ with the highest root of $C$ is in
$R_1(v)$ while the other points of $C$ are in $R_0(v)$.

Proof. Let $u$ be as in (1) and put $r=u-v$. The vector $r$ is
orthogonal to $u$ and has inner product $\le 0$ with all roots of
$R_0$ (because $-v$ does) so it is a highest root of some component $C$
of $R_0(u)$. The vector $r$ therefore determines some simple spherical
Dynkin diagram $C$ contained in $\Lambda$. Any root $c$ of $C$ has
$(c,v+r)=(c,u)=0$, so $c$ is in $R_i(v)$ where $i=(c,r)$. This gives a
map from (1) to (2).

Conversely if we start with a Dynkin diagram $C$ satisfying (2) and
put $u=v+r$ (where $r$ is the highest root of $C$) then $(c,u)=0$ for
all $c$ in $C$, so $(r,u)=0$ as $r$ is a sum of the $c$'s. This
implies that $u^2=2n$ and $(u,v)=2n$. We now have to show that $u$
is in $D$. Let $s$ be any simple root of $D$. If $s$ is in $C$ then
$(s,r)=-(s,v)$ and if $s$ is not in $C$ then $(s,r)\ge 0$, so in any
case $(s,u)=(s,v+r)\ge 0$ and hence $u$ is in $D$. This gives a map
from (2) to (1) and shows that (1) and (2) are equivalent.

Condition (3) is just the condition (2) written out explicitly for
each possible $C$, so (2) and (3) are also equivalent. This proves
theorem 2.2.

We  define the height of a vector $u$ in
$II_{1,25}$ to be $(u,w)$. We show how to calculate
the heights of vectors of $II_{1,25}$ that have been found with the
algorithm above.

\proclaim Lemma 2.3.
Suppose $u$, $v$ are vectors in $D$ of norms $2n$, $2(n-1)$ with
$(u,v)=2n$ and suppose that $v=u-r$ for some root $r$ of $u^\perp$
corresponds to the component $C$ of $R_0(u)$. Then
$$\height(u)=\height(v)+h-1$$
where $h$ is the Coxeter number of the component $C$. 

Proof. We have $v=u-r$ where $r$ is the highest root of $C$, so
$\height(u)=\height(v)+(r,w)$. We have $r=\sum_im_ic_i$ where the
$c_i$ are the simple roots of $C$ with weights $m_i$ and
$\sum_im_i=h-1$. All the $c_i$ have inner product $1$ with $w$, so
$(r,w)=h-1$. This proves lemma 2.3.

\proclaim Lemma 2.4.
Let $u$ be a primitive vector of $D$ such that there
is a norm 0 vector $z$ with $(z,u)=0$ or $1$, and suppose that $z$
corresponds to a Niemeier lattice $B$ with Coxeter number $h$.
\item{(1)} 
If $u$ has norm 0 then its height is $h$. The Dynkin diagram of
$u^\perp$ is the extended Dynkin diagram of $B$.
\item{(2)}
If $u$ has positive norm  then
$\height(u)=1+(1+u^2/2)h$. The Dynkin diagram of $u^\perp$ is the
Dynkin diagram of $B$ if $u^2>2$ and the Dynkin diagram of $B$ plus an
$a_1$ if $u^2=2$.

Proof.
\item{(1)}
The Dynkin diagram of $u^\perp$ is a union of extended Dynkin
diagrams. If this union is empty then $u$ must be $w$ and therefore
has height $0=h$. If not then let $C$ be one of the components. We
have $u=\sum_im_ic_i$ where the $c_i$'s are the simple roots of $C$
with weights $m_i$. Also $\sum_im_i=h$ because $C$ is an extended
Dynkin diagram and all the $c_i$'s have height 1, so $u$ has height
$h$.
\item{(2)}
As $u$ has inner product 1 with a norm 0 vector $z$ of $D$ we can put
$u=nz+z'$ with $u^2=2n$ and $z'^2=0$, $(z,z')=1$.  By part (1) $z$ has
height $h$. We have $z'=z+r$ where $r$ is a simple root of $D$, so
$\height(z')=\height(z)+\height(r)=h+1$. Hence
$\height(u)=nh+h+1=1+(1+u^2/2)h$. The lattice $u^\perp$ is $B\oplus N$
where $N$ is a one dimensional lattice of determinant $2n$, so the
Dynkin diagram is that of $B$ plus that of $N$, and the Dynkin diagram
of (norm 2 roots of) $N$ is empty unless $2n=2$ in which case it is
$a_1$. This proves that the Dynkin diagram of $u^\perp$ is what it is
stated to be. This proves lemma 2.4.

Orbits of norm $2$ vectors $u\in II_{1,25}$ correspond to even 25
dimensional positive definite lattices $B$ of determinant 2, where
$B(-1)\cong u^\perp$.  One part of the algorithm for finding vectors
of norm $2n$ consists of finding the vectors $u$ such that there are
no roots in $u^\perp$.  For norm $2$ vectors $u$ the following lemma shows
that there are no such vectors.

\proclaim Lemma 2.5.
If $u\in II_{1,25}$ has norm 2 then $u^\perp$ contains roots. In other
words every 25 dimensional even positive definite lattice of
determinant 2 has a root.

Proof. If $u^\perp$ contains no roots then, by the algorithm of
section 1, $u=w+u_1$ for some $u_1$ in $D$. We have
$u_1^2=u^2-2\height (u)$, so $u_1^2=0$ and $u$ has height 1 because
$u_1^2\ge 0$, $u^2=2$ and the height of $u$ is positive. Then
$\height(u_1)=\height(u)=1$, so $u$ is a norm 0 vector in $D$ that has
inner product $1$ with the norm 0 vector $w$ of $D$, but this is
impossible as $u-w$ would be a norm $-2$ vector separating the two
vectors $u$ and $w$ of $D$. This proves lemma 2.5.

\proclaim Theorem 2.6. 
Suppose that $u\in D$ has norm 2.  Then $$w=\rho+\height(u)u/2$$ where
$\rho$ is the Weyl vector of the root system of $u^\perp$.  Also
$-2\rho^2=\height(u)^2$.

Proof. The vector $w$ is fixed by any automorphism fixing $D$, so by
lemma 1.4 the vector $w$ must be in the space spanned by $u$ and the
roots of $u^\perp$.  However $w$ also has inner product $1$ with all
simple roots of $u^\perp$ and has inner product $\height(u)$ with $u$,
so $w$ must be $\rho+\height(u)u/2$.  Taking norms of both sides of
$w=\rho+\height(u)u/2$, and using the facts that $w^2=0$,
$(u,\rho)=0$, and $(u,u)=2$, shows that $-2\rho^2=\height(u)^2$.  This
proves theorem 2.6.

In particular we find the strange consequence that the norm of the
Weyl vector of any 25 dimensional even positive definite lattice of
determinant 2 must be a half a square.

Norm $4$ vectors in the fundamental domain $D$ of $II_{1,25}$
correspond to 25 dimensional unimodular lattices $A=A_1\oplus I^n$,
where $u^\perp$ is the lattice of even elements of $A(-1)$ and
$A_1$ has no norm 1 vectors . The odd
vectors of $A(-1)$ can be taken as the projections of the vectors $y$
with $(y,u)=2$ into $u^\perp$.  A norm 4 vector $u$ can behave in 4
different ways, depending on whether the unimodular lattice $A_1$ with
no norm 1 vectors corresponding to $u$ is at most 23 dimensional, or
24 dimensional and odd, or 24 dimensional and even, or 25 dimensional.

\proclaim Theorem 2.7.
Norm 1 vectors of $A$ correspond to norm 0 vectors $z$ of $II_{1,25}$
with $(z,u)=2$. Write $A=A_1\oplus I^n$ where $A_1$ has no vectors of
norm 1. Then $u$ is in exactly one of the following four classes:
\item{(1)} 
$u$ has inner product 1 with a norm 0 vector. The lattice $A_1$ is a
Niemeier lattice.
\item{(2)}
$A$ has at least 4 vectors of norm 1, so that $A_1$ is at most 23
dimensional (but may be even). There is a unique norm 0 vector $z$ of
$D$ with $(z,u)=2$ and this vector $z$ is of the same type as either
of the two even neighbors of $A_1\oplus I^{n-1}$.
\item{(3)}
$A_1$ is 24 dimensional and odd. There are exactly two norm 0 vectors
that have inner product $2$ with $u$, and they are both in $D$. They
have the types of the two even neighbors of $A_1$.
\item{(4)}
$A=A_1$ has no vectors of norm 1.

Proof. The vector $z$ is a norm 0 vector with $(z,u)=2$ if and only if
$u/2-z$ is a norm 1 vector of $A$. Most of 2.7 follows from this. The
only non-trivial things to check are the statements about norm 0
vectors that are in $D$.

If $u$ does not have inner product 1 with any norm 0 vector then a
norm 0 vector $z$ with $(z,u)=2$ is in $D$ if and only if it has inner
product $\ge 0$ with all simple roots of $u^\perp$, so there is one
such vector in $D$ for each orbit of such norm 0 vectors under the
reflection group of $u^\perp$. If $A$ has at least 4 vectors of norm 1
then they form a single orbit under the Weyl group of (the norm 2
vectors of) $u^\perp$, which proves (2), while if $A$ has only two
vectors of norm 1 then they are both orthogonal to all norm 2 vectors
of $A$ and so form two orbits under they Weyl group of $u^\perp$. This
proves theorem 2.7.

\proclaim Theorem 2.8. 
Suppose that $u$ is a norm 4 vector corresponding to a unimodular 25
dimensional lattice $A=A_1\oplus I^{25-n}$ with $2n\ge 4$ vectors of
norm 1. Let $\rho$ be the Weyl vector of the root system of norm $-2$
roots of $u^\perp$ (which is the Weyl vector of the norm $-2$ vectors
of $A(-1)$) and let let $h$ be the Coxeter number of the even
neighbors of the 24 dimensional unimodular lattice $A_1\oplus
I^{24-n}$. Then $\height(u)=(w,u)=2(h+n-1)$, $w=\rho+\height(u)u/4$,
and $-\rho^2=(h+n-1)^2$.

Proof.  There is a unique norm 2 vector $z$ of $D$ with $(z,u)=2$; we
let $i$ be its projection into $u^\perp$.  The lattice $A$ has at
least 4 vectors of norm 1, so any vector of norm 1 and in particular
$i$ is in the vector space generated by vectors of norm $-2$ of
$u^\perp$. Hence by lemma 1.5 and the same argument as in theorem 2.6
we have $w=\rho+\height(u)u/4$. The norm $-4$ vector $2i$ of $u^\perp$ is
the sum of $-2(n-1)$ simple roots of the $d_n$ component of the Dynkin
diagram of $u^\perp$, so $(2i,w)=(2i,\rho)=-2(n-1)$.

The vector $i$ is the projection of $z$ into $u^\perp$, so $i=z-u/2$,
and hence
$$\eqalign{
\height(u)&=(w,u)\cr
&= 2(w,z-i)\cr
&= 2(\height(z)+n-1)\cr
&=2(h+n-1).\cr
}$$
If we calculate the norms of both sides of $w=\rho+\height(u)u/4$ we
find that $-\rho^2=(h+n-1)^2$. This proves theorem 2.8.

\noindent
{\bf Example 2.9.} Suppose $u$ corresponds to the lattice $I^{25}$.
The number $n$ is then 25 and the root system of the norm 2 vectors is
$D_{25}$, so the Weyl vector $\rho$ can be taken as
$(0,1,2,\ldots,24)$. The even neighbors of $I^{24}$ are both $D_{24}$
with Coxeter number $h=46$, so we find that
$0^2+1^2+2^2+\cdots+24^2=\rho^2=(h+n-1)^2=70^2$. Watson showed that
the only solution of $0^2+1^2+\cdots+k^2=m^2$ with $k\ge 2$ is $k=24$.
See [C-S Chapter 26] for a construction of the Leech lattice using this 
equality. 

\proclaim Theorem 2.10. 
Suppose that $u$ is a norm $4$ vector of $D$ with exactly two norm 0
vectors $z_1$, $z_2$ that have inner product $2$ with $u$, and suppose
that there are no norm 0 vectors that have inner product 1 with
$u$. Then $z_1$ and $z_2$ are both in $D$ and have Coxeter numbers
$h_1$, $h_2$ where $h_i=(z_i,w)$.  Then
$$w=\rho+(h_1z_2+h_2z_1)/2$$
where $\rho$ is the Weyl vector of the norm $-2$ vectors of $u^\perp$. 
Also $u=z_1+z_2$, $\height(u)=h_1+h_2$,
$-\rho^2=h_1h_2$, and $u^\perp$ has $8(h_1+h_2-2)$ roots.

Proof. The vector $u-z_1$ is a norm 0 vector which has inner product
$2$ with $u$ and so must be $z_2$. Hence $u=z_1+z_2$ and
$\height(u)=\height(z_1)+\height(z_2)=h_1+h_2$.

There is a norm 0 vector that has inner product $2$ with $u$, and any
automorphism of $L$ fixing $D$ also fixes $w$, so by lemma 1.5 $w$ is
a linear combination of $z_1$, $z_2$, and the roots of $R$ in
$u^\perp$. Using the facts that $(w,z_1)=h_1$, $(w,z_2)=h_2$, and
$(w,r)=-r^2/2$ for any simple root $r$ in $u^\perp$ shows that $w$
must then be $\rho+(h_1z_2+h_2z_1)/2$. Using the fact that $w^2=0$
this shows immediately that $-\rho^2=h_1h_2$. The number of roots
follows from remark 2.12 below.  This proves theorem 2.10.

\proclaim Corollary 2.11.
If $A_1$ is an odd 24 dimensional positive definite unimodular lattice
with no vectors of norm 1 and whose even neighbors have Coxeter
numbers $h_1$ and $h_2$, then $\rho^2=h_1h_2$ where $\rho$ is the Weyl
vector of $A_1$.

Proof. This follows immediately from theorem 2.10, using the fact that
$A_1\oplus I$ is the 25 dimensional unimodular lattice corresponding to
$u$ as in 2.10.

Remark.  Let $B_1$, $B_2$ be the two even neighbors of $A_1$. Then it is
not hard to show that $h_2\le 2h_1+2$, and there are several lattices
$A_1$ for which equality holds.

\noindent
{\bf Remark 2.12.}  Theorem 13.1 and corollary 13.2 of [B95] show that
the height of a vector in the fundamental domain of $II_{1,25}$ can be
written as an explicit linear combination of the theta functions of
cosets of the lattice $u^\perp$.  In particular we find that if $u$ is
a norm 2 vector then
$$12\height(u)=18-4z_1+r$$ where $r$ is the number of norm $-2$ vectors
of $u^\perp$ and $z_i$ is the number of norm 0 vectors having inner
product $i$ with $u$ (so $z_1$ is 0 or 2 and is 2 if and only if the
lattice $u^\perp$ is the sum of a one dimensional lattice and an even
lattice). Similarly if $u$ has norm $4$ and corresponds to a 25
dimensional unimodular lattice $A$ then
$$8t=20-2z_2-8z_1+r$$ where $r$ is the number of norm 2 vectors of
$A$, $z_2$ is the number of norm 1 vectors of $A$, and $z_1$ is 1 if
$A$ is the sum of a Niemeier lattice and a one dimensional lattice and
is 0 otherwise.  Note that these relations give congruences for the
numbers of roots that immediately imply that 25 dimensional even
lattices of determinant 2 and 25 dimensional unimodular lattices
always have roots. There are similar relations and congruences for
larger norm vectors of $II_{1,25}$.

There are several other genuses of lattices that can be classified
using $II_{1,25}$. Most of these do not seem important enough to be
worth publishing, but here is a summary of what is available just in
case anyone finds a use for any of these. The 24 dimensional even
positive definite lattices of determinant 5 are easy to classify as
they turn out to correspond to pairs consisting of a norm $2$ vector
$u$ of $II_{1,25}$ together with a norm $-2$ root $r$ with $(r,u)=1$,
and these can easily be read off from the list of norm 2 vectors.  The
25 dimensional positive definite even lattices of determinant 6
correspond to the norm 6 vectors in $II_{1,25}$ and can be classified
from the norm 4 vectors using the algorithm; there are 2825 orbits if
I have made no mistakes. A list of them is available from my home
page. These can be used to classify the 26 dimensional even positive
definite lattices of determinant 3, because the norm 2 roots of such
lattices correspond to the norm 6 vectors of $II_{1,25}$. (There is a
unique such lattice with no roots; see the next section.) There are
between 677 and 681 such lattices, and a provisional list is available
from my home page (there are a few small ambiguities that I have not
yet got around to resolving). If such a lattice has no norm 6 roots
then the number of norm 2 vectors is divisible by 6.  With a lot more
effort it should be possible to classify the 26 dimensional unimodular
lattices by finding the (roughly 50000?) orbits of norm 10 vectors of
$II_{1,25}$; see the next section.

\proclaim 3.~Lattices with no roots.

In this section we show that there is a unique 26 dimensional positive
definite unimodular lattice with no roots. Conway and Sloane use this
result in their proof [C-S98] that there is a positive definite
unimodular lattice with no roots in all dimensions greater than 25. We
also show that the number of norm 2 vectors of a 26 dimensional
unimodular lattice is divisible by 4, and sketch a construction of a
27 dimensional unimodular lattice with no roots.

\proclaim Lemma 3.1. 
A 26-dimensional unimodular lattice $L$ with no vectors of norm 1 has
a characteristic vector of norm 10.

Proof. If $L$ has a characteristic vector $x$ of norm 2 then $x^\perp$
is a 25 dimensional even lattice of determinant 2 and therefore has a
root $r$ by theorem 2.6; $2r+x$ is a characteristic vector of norm
10. If the lemma is not true we can therefore assume that $L$ has no
vectors of norm 1 and no characteristic vectors of norm 2 or 10. Its
theta function is determined by these conditions and turns out to be
$1-156q^2+\cdots$ which is impossible as the coefficient of $q^2$ is
negative. This proves lemma 3.1.

\proclaim Lemma 3.2. There is a bijection between isomorphism classes of
\item{(1)} 
Norm 10 characteristic vectors $c$ in 26-dimensional positive definite
unimodular lattices $L$, and
\item{(2)}
Norm $10$ vectors $u$ in $II_{1,25}$
\item{}
given by $c^\perp(-1)\cong u^\perp$. 
\item{}We have $\Aut(L,c)=\Aut(II_{1,25}, u)$. 

Proof. Routine. Note that $-1$ is a square mod 10. This proves lemma 3.2. 

Lemmas 3.1 and 3.2 give an algorithm for finding 26 dimensional
unimodular lattices $L$. It is probably not hard to implement this on
a computer if one is given a computer algorithm for deciding when 2
vectors of the Leech lattice are conjugate under its automorphism
group; such an algorithm has been described by Allcock in [A]. The
main remaining open problem is to find a use for these lattices! We
now apply this algorithm to find the unique such lattice with no
roots.

\proclaim Lemma 3.3. 
Take notation as in lemma 3.2. The lattice $L$ has no roots if and
only if $u^\perp$ has no roots and $u$ does not have inner product 1,
2, 3, or 4 with any norm 0 vector.

Proof. If $u^\perp$ has roots then obviously $L$ has too. If there is
a norm 0 vector $z$ that has inner product 1, 2, 3, or 4 with $u$ then
the projection $z_u$ of $z$ into $u^\perp$ has norm $-1/10$, $-4/10$,
$-9/10$, or $-16/10$. The lattice $L(-1)$ contains $u^\perp+c$, and
the vector $z_u\pm 3c/10$, $z_u\pm 4c/10$, $z_u\pm c/10$, or $z_u\pm
2c/10$ is in $L$ for some choice of sign and has norm $-1$, $-2$,
$-1$, or $-2$. Hence if $u$ has inner product 1, 2, 3, or 4 with some
norm 0 vector then $L$ has roots. Conversely if $L$ has a root $r$
then either $r$ has norm 2 and inner product 0, $\pm 2$, $\pm 4$ with
$c$ or it has norm 1 and inner product $\pm 1$, $\pm 3$ with $c$, and
each of these cases implies that $u^\perp$ has roots or that $u$ has
inner product 1, 2, 3, or 4 with some norm 0 vector by reversing the
argument above. This proves lemma 3.3.

Now let $L$ be a 26 dimensional unimodular lattice with no roots
containing a characteristic vector $c$ of norm 10, and let $u$ be a
norm $10$ vector of $D$ corresponding to it as in 3.2.

\proclaim Lemma 3.4. 
$u=z+w$, where $z$ is a norm 0 vector of $D$ corresponding to a
Niemeier lattice with root system $A_4^6$, and $w$ is the Weyl vector
of $D$. In particular $u$ is determined up to conjugacy under
$\Aut(D)$.

Proof. The lattice $u^\perp$ has no roots so $u=w+z$ for some vector
$z$ of $D$. By lemma 3.3 $u$ does not have inner product 1, 2, 3, or 4
with any norm 0 vector, so $(z,w)=(u,w)\ge 5$. Hence
$$10=u^2=z^2+2(z,w)\ge 2(z,w)\ge 10$$ so $(z,w)=5$ and $z^2=0$. The
only norm 0 vectors $z$ in $D$ with $(z,w)=5$ are the primitive ones
corresponding to $A_4^6$ Niemeier lattices, which form one orbit under
$\Aut(D)$. This proves lemma 3.4.

\proclaim Lemma 3.5. 
If $u=z+w$ is as in lemma 3.4 then the 26 dimensional unimodular lattice
corresponding to $u$ has no roots.

Proof. The lattice $u^\perp$ obviously has no roots so by lemma 3.3 we
have to check that there are no norm 0 vectors that have inner product
$1$, $2$, $3$, or $4$ with $u$. Let $x$ be any norm 0 vector in the
positive cone. If $x$ has type $A_4^6$ then $(x,u)\ge (x,w)\ge 5$; if
$x$ has Leech type then $(x,u)\ge (x,z)\ge 5$; if $x$ has type
$A_1^{24}$ then $(x,u)=(x,w)+(x,z)\ge 2+3=5$ ($(x,z)$ cannot be $2$ as
there are no pairs of norm 0 vectors of types $A_1^{24}$ and $A_4^6$
that have inner product $2$ by the classification of 24 dimensional
unimodular lattices); and if $x$ has any other type then
$(x,u)=(x,w)+(x,z)\ge 3+2=5$. This proves lemma 3.5.

\proclaim Theorem 3.6. 
There is a unique 26 dimensional positive definite unimodular lattice
$L$ with no roots. Its automorphism group is isomorphic to the group
$O_5(5)=2.G.2$ of order $2^8.3^2.5^4.13$ and acts transitively on the
624 characteristic norm 10 vectors of $L$.

Proof. By lemma 3.1 $L$ has a characteristic vector of norm 10, so by
lemmas 3.3 and 3.4 $L$ is unique and its automorphism group acts
transitively on the characteristic vectors of norm 10. By lemma 3.5
$L$ exists. The theta function is determined by the conditions that
$L$ has no vectors of norm 1 or 2 and no characteristic vectors of
norm 2, and it turns out that the number of characteristic vectors of
norm 10 is 624. The stabilizer of such a vector is isomorphic to
$\Aut(II_{1,25}, u)$, which is a group of the form $5^3.2.S_5$ where
$S_5$ is the symmetric group on 5 letters.  This determines the order
of the automorphism group of the lattice.  From this it is not
difficult to determine it precisely; we omit the details.  This proves
theorem 3.6.

We now show that the number of norm 2 vectors of any 26 dimensional
even positive definite unimodular lattice is divisible by 4.  There
are strictly 26 dimensional unimodular lattices with no roots or with
4 roots, so this is the best possible congruence. For unimodular
lattices of dimension less than 26 there are congruences modulo higher
powers of 2 for the number of roots.

\proclaim Lemma 3.7. 
If $L$ is a 25-dimensional positive definite lattice of determinant 2
then the number of norm 2 roots of $L$ is $2\bmod 4$.

Proof. The even vectors of $L$ form a lattice isomorphic to the
vectors that have even inner product with some vector $b$ in an even
25-dimensional lattice $B$ of determinant 2. (Note that $b$ is not in
$B'-B$.) The number of roots of $B$ is $12t-10$ or $12t-18$ where $t$
is the height of the norm $2$ vector of $D$ corresponding to $B$ by
remark 2.12, so it is sufficient to prove that the number of norm 2
vectors of $B$ that have odd inner product with $b$ is divisible by 4.

The vector $b$ has zero inner product with $u$ and integral inner
product with $w$, so by theorem 2.6 $b$ has integral inner product
with $\rho$. Hence $b$ has even inner product with the sum of the
positive roots of $B$, so it has odd inner product with an even number
of positive roots. This implies that the number of roots of $B$ that
have odd inner product with $b$ is divisible by 4. This proves lemma
3.7.

\proclaim Corollary 3.8. 
If $L$ is a 26 dimensional unimodular lattice then the number of norm
2 vectors of $L$ is divisible by 4.

Proof. The result is obvious if $L$ has no norm 2 roots, so let $r$ be
a norm 2 vector of $L$. The lattice $r^\perp$ is a 25 dimensional even
lattice of determinant so by remark 2.12 the number of roots of
$r^\perp$ is $2\bmod 4$. The number of roots of $L$ not in $r^\perp$
is $4h-6$ where $h$ is the Coxeter number of the component of $L$
containing $r$, so the number of norm 2 vectors of $L$ is divisible by
4. This proves corollary 3.8.

{\bf Remark.}  A similar but more complicated argument can be used to
show that there is a unique even 26 dimensional positive definite
lattice of determinant 3 with no roots. Gluing on a one dimensional
lattice to this gives a unique 27 dimensional unimodular lattice with no
roots and a characteristic vector of norm 3. As a different proof of
this has already been published in [E-Z] we will just give a brief
sketch of the proof from [B]. (The preprint [B-V] shows that there are
exactly three 27 dimensional positive definite unimodular lattices with no
roots.)  Let $L$ be a 27 dimensional positive definite unimodular
lattice with no roots and a characteristic vector $c$ of norm 3. The
theta function of $L$ is determined by these conditions and this
implies that $L$ has vectors of norm 5; let $v$ be such a vector.
Then $\langle v,c\rangle^\perp$ is a 25 dimensional even lattice $X$
of determinant 14 such that $X'/X$ is generated by an element of norm
$1/14\bmod 2$. Such lattices $X$ correspond to norm $14$ vectors $x$
in the fundamental domain $D$ of $II_{1,25}$, and the condition that
$L$ has no vectors of norm 1 or 2 implies that there are exactly two
possibilities for $x$: $x$ is either the sum of $w$ and a norm 0
vector of height 7 corresponding to $A_6^4$, or $x$ is the sum of $w$
and a norm $2$ vector of height 6 corresponding to the 25 dimensional
lattice of determinant 2 with root system $a_2^9$. Both of these $x$'s
turn out to give the same lattice $L$, which therefore has two orbits
of norm 5 vectors and is the unique 27 dimensional positive definite
unimodular lattice with no roots and with characteristic vectors of
norm 3.

\proclaim Table 0.~The primitive norm $0$  vectors of $II_{1,25}$. 

We list the set of orbits of primitive norm 0 vectors $z$ of
$II_{1,25}$, which is of course more or less the same as the well
known list of Niemeier lattices (see [C-S table 16.1]).  The height is
just $(w,z)$ where $w$ is the Weyl vector of a fundamental domain
containing $z$. The letter after the height is just a name to
distinguish vectors of the same height, and is the letter referred to
in the column headed ``Norm $0$ vectors'' of table $1$.  The column
headed ``Group'' is the order of the subgroup of $\Aut(D)$ fixing the
primitive norm 0 vector. However note that the group order is {\it
not} (usually) the order of the quotient of the automorphism group of
the Niemeier lattice by the reflection group; see [C-S chapter 16] for
a description of the relation between these groups. For the vector $w$
of height 0 the group is the infinite group of automorphisms of the
affine Leech lattice and is an extension of a finite group of the
order given by the group of translations of the Leech lattice
$\Lambda$.

\bigskip
\halign{
\hfill#&#&~~\hfill$#$&~~\hfill#\cr
Height&&\hbox{Roots}&Group\cr
\cr
0&x&\hbox{None}&$\Lambda\cdot$8315553613086720000\cr
2&a&A_1^{24}&1002795171840\cr
3&a&A_2^{12}&138568320\cr
4&a&A_3^8   & 688128\cr
5&a&A_4^6& 30000\cr
6&d&D_4^6&138240\cr
6&a&A_5^4D_4& 3456\cr
7&a&A_6^4&  1176\cr
8&a&A_7^2D_5^2&256\cr
9&a&A_8^3&324\cr
10&d&D_6^4&384\cr
10&a&A_9^2D_6&80\cr
12&e&E_6^4&432\cr
12&a&A_{11}D_7E_6&24\cr
13&a&A_{12}^2&52\cr
14&d&D_8^3&48\cr
16&a&A_{15}D_9&16\cr
18&d&D_{10}E_7^2&8\cr
18&a&A_{17}E_7&12\cr
22&d&D_{12}^2&8\cr
25&a&A_{24}&10\cr
30&e&E_8^3&6\cr
30&d&D_{16}E_8&2\cr
46&d&D_{24}&2\cr
}

\proclaim Table $1$.~The norm $2$  vectors of $II_{1,25}$. 

The following sets are in natural 1:1 correspondence:
\item{(1)} 
Orbits of norm $2$ vectors in $II_{1,25}$ under $\Aut(II_{1,25})$.
\item{(2)} 
Orbits of norm $2$ vectors $u$ of $D$ under $\Aut(D)$. 
\item{(3)} 
25 dimensional even bimodular lattices $L$. 

The lattice $L(-1)$ is isomorphic to $u^\perp$. Table 1 lists the
121 elements of any of these three sets.

The {\it height}  is the height of the norm $2$ vector $u$ of $D$,
in other words $(u,w)$ where $w$ is the Weyl vector of $D$. The
letter after the height is just a name to distinguish vectors of the
same height, and is the letter referred to in the column headed ``Norm
$2$'' of table $2$. An asterisk after the letter means that the
vector $u$ is of type 1, in other words the lattice $L$ is the sum of
a Niemeier lattice and $a_1$.

The column ``Roots'' gives the Dynkin diagram of the norm 2 vectors of
$L$ arranged into orbits under $\Aut(L)$. ``Group'' is the order of
the subgroup of $\Aut(D)$ fixing $u$. The group $\Aut(L)$ is a split
extension $R.G$ where $R$ is the Weyl group of the Dynkin diagram and
$G$ is isomorphic to the subgroup of $\Aut(D)$ fixing $u$.

``$S$'' is the maximal number of pairwise orthogonal roots of $L$. 

The column headed ``Norm 0 vectors'' describes the norm 0 vectors $z$
corresponding to each orbit of roots of $u^\perp$ where $u$ is in $D$.
A capital letter indicates that the corresponding norm 0 vector is
twice a primitive vector, otherwise the norm 0 vector is
primitive. $x$ stands for a norm 0 vector of type the Leech
lattice. Otherwise the letter $a$, $d$, or $e$ is the first letter of
the Dynkin diagram of the norm 0 vector, and its height is given by
${\rm height}(u)-h+1$ where $h$ is the Coxeter number of the component
of the Dynkin diagram of $u$.

For example, the norm $2$ vector of type $23a$ has 3 components in
its root system, of Coxeter numbers 12, 12, and 6, and the letters are
$e$, $a$, and $d$, so the corresponding norm 0 vectors have Coxeter
numbers 12, 12, and 18 and hence are norm 0 vectors with Dynkin
diagrams $E_6^4$, $A_{11}D_7E_6$, and $D_{10}E_7^2$.

\bigskip
\halign{
\hfill#&#\hfill&~\hfill$#$&~~\hfill#&~~~\hfill#&\hfill#&#\cr
Height&&\hbox{Roots}&Group&$S$&~Norm 0&~vectors\cr
\cr
1&a*&a_1 &8315553613086720000 &1&X\cr
\cr
2&a &a_2 &991533312000 &1&x\cr
\cr
3&a &a_1^9 &92897280 &9&a\cr
\cr
4&a &a_2a_1^{12} &190080 &13&aa\cr
\cr
5&a*&a_1^{24}a_1 &244823040 &25&aA\cr
5&b &a_2^4a_1^9 &3456 &13&aa\cr
5&c &a_3a_1^{15} &40320 &17&aa\cr
\cr
6&a &a_2^9 &3024 &9&a\cr
6&b &a_3a_2^5a_1^6 &240 &13&aaa\cr
\cr
7&a*&a_2^{12}a_1 &190080 &13&aA\cr
7&b &a_3^3a_2^4a_1^3 &48 &13&aaa\cr
7&c &a_3^4a_1^8a_1 &384 &17&aad\cr
7&d &a_4a_2^6a_1^5 &240 &13&aaa\cr
7&e &d_4a_1^{21} &120960 &25&ad\cr
\cr
8&a &a_3^6a_2 &240 &13&ad\cr
8&b &a_4a_3^3a_2^3a_1^2 &12 &13&aaaa\cr
8&c &d_4a_2^9 &864 &13&aa\cr
\cr
9&a*&a_3^8a_1 &2688 &17&aA\cr
9&b &a_4^2a_3^4a_1 &16 &13&aaa\cr
9&c &a_4^3a_3a_2^2a_1^3 &12 &13&aaaa\cr
9&d &d_4a_3^4a_3a_1^3 &48 &17&aada\cr
9&e &a_5a_3^3a_2^4 &24 &13&aaa\cr
9&f &a_5a_3^4a_1^6 &48 &17&aaa\cr
\cr
10&a &d_4a_4^3a_2^3 &12 &13&aaa\cr
10&b &a_5a_4^2a_3^2a_2a_1 &4 &13&aaaaa\cr
\cr
11&a*&a_4^6a_1 &240 &13&aA\cr
11&b &d_4^4a_1^9 &432 &25&dd\cr
11&c &a_5d_4^2a_3^3 &24 &17&daa\cr
11&d &a_5a_5a_4^2a_3a_1 &4 &13&aaaaa\cr
11&e &a_5^2d_4a_3^2a_1^2a_1 &8 &17&aaaad\cr
11&f &a_5^3a_2^4 &48 &13&aa\cr
11&g &d_5a_3^6a_1 &48 &17&aad\cr
11&h &a_6a_4^2a_3^2a_2a_1 &4 &13&aaaaa\cr
\cr
12&a &a_5^4a_2 &24 &13&ad\cr
12&b &d_5a_4^4a_2 &8 &13&aaa\cr
12&c &a_6d_4a_4^3 &6 &13&aaa\cr
12&d &a_6a_5^2a_3a_2^2 &4 &13&aaaa\cr
\cr
13&a*&a_5^4d_4a_1 &48 &17&aaA\cr
13&b &d_5a_5^2d_4a_3a_1 &4 &17&aaada\cr
13&c &d_5a_5^3a_1^3a_1 &12 &17&aaae\cr
13&d*&d_4^6a_1 &2160 &25&dD\cr
13&e &a_6^2a_5a_4a_1^2 &4 &13&aaaa\cr
13&f &a_7a_5a_4^2a_3 &4 &13&aaaa\cr
13&g &a_7a_5d_4a_3^2a_1^2 &4 &17&aaaaa\cr
\cr
14&a &a_6a_6d_5a_4a_2 &2 &13&aaaaa\cr
14&b &a_6^3d_4 &12 &13&aa\cr
14&c &a_7a_6a_5a_4a_1 &2 &13&aaaaa\cr
\cr
15&a*&a_6^4a_1 &24 &13&aA\cr
15&b &d_5^3a_5a_3 &12 &17&ade\cr
15&c &d_6d_4^4a_1^3 &24 &25&ddd\cr
15&d &d_6a_5^2a_5a_3 &4 &17&aada\cr
15&e &a_7d_5^2a_3^2a_1 &4 &17&aaad\cr
15&f &a_7^2d_4^2a_1 &8 &17&aad\cr
15&g &a_8a_5^3 &6 &13&aa\cr
15&h &a_8a_6a_5a_3a_2 &2 &13&aaaaa\cr
\cr
16&a &a_7^3a_2 &12 &13&ad\cr
16&b &a_8a_6d_5a_4 &2 &13&aaaa\cr
\cr
17&a*&a_7^2d_5^2a_1 &8 &17&aaA\cr
17&b &e_6a_5^3d_4 &12 &17&aae\cr
17&c &a_7d_6d_5a_5 &2 &17&daaa\cr
17&d &a_7^2d_6a_3a_1 &4 &17&aada\cr
17&e &a_8a_7^2a_1 &4 &13&aaa\cr
17&f &a_9d_5a_5d_4a_1 &2 &17&aaaaa\cr
17&g &a_9a_7a_4^2 &4 &13&aaa\cr
\cr
18&a &e_6a_6^3 &6 &13&aa\cr
18&b &a_9a_8a_5a_2 &2 &13&aaaa\cr
\cr
19&a*&a_8^3a_1 &12 &13&aA\cr
19&b &d_6^3d_4a_1^3 &6 &25&ddd\cr
19&c &a_7e_6d_5^2a_1 &4 &17&eaad\cr
19&d &d_7a_7d_5a_5 &2 &17&aaad\cr
19&e &d_7a_7^2a_3a_1 &4 &17&aaad\cr
19&f &a_9a_7d_6a_1a_1 &2 &17&aaaad\cr
19&g &a_{10}a_7a_6a_1 &2 &13&aaaa\cr
\cr
20&a &a_8^2e_6a_2 &4 &13&aaa\cr
20&b &a_{10}a_8d_5 &2 &13&aaa\cr
\cr
21&a*&a_9^2d_6a_1 &4 &17&aaA\cr
21&b &a_{11}d_6a_5a_3 &2 &17&aaaa\cr
21&c &a_{11}a_8a_5 &2 &13&aaa\cr
21&d*&d_6^4a_1 &24 &25&dD\cr
21&e &a_9e_6d_6a_3 &2 &17&aaad\cr
\cr
23&a &d_7e_6^2a_5 &4 &17&ead\cr
23&b &d_8d_6^2d_4a_1 &2 &25&dddd\cr
23&c &a_9d_7^2 &4 &17&da\cr
23&d &a_9d_8a_7 &2 &17&daa\cr
23&e &a_{11}d_7d_5a_1 &2 &17&aaad\cr
\cr
24&a &a_{11}^2a_2 &4 &13&ad\cr
24&b &a_{12}e_6a_6 &2 &13&aaa\cr
\cr
25&a*&a_{11}d_7e_6a_1 &2 &17&aaaA\cr
25&b &a_{13}d_6d_5 &2 &17&aaa\cr
25&e*&e_6^4a_1 &48 &17&eE\cr
\cr
26&a &a_{13}a_{10}a_1 &2 &13&aaa\cr
\cr
27&a*&a_{12}^2a_1 &4 &13&aA\cr
27&b &e_7d_6^3 &3 &25&dd\cr
27&c &a_9a_9e_7 &2 &17&ada\cr
27&d &d_9a_9e_6 &2 &17&ada\cr
27&e &a_{11}d_9a_5 &2 &17&aad\cr
27&f &a_{14}a_9a_2 &2 &13&aaa\cr
\cr
29&a &a_{11}e_7e_6 &2 &17&daa\cr
29&d*&d_8^3a_1 &6 &25&dD\cr
\cr
31&a &d_8^2e_7a_1a_1 &2 &25&ddde\cr
31&b &d_{10}d_8d_6a_1 &1 &25&dddd\cr
31&c &a_{15}d_8a_1 &2 &17&aad\cr
\cr
33&a*&a_{15}d_9a_1 &2 &17&aaA\cr
33&b &a_{15}e_7a_3 &2 &17&aad\cr
33&c &a_{17}a_8 &2 &13&aa\cr
\cr
35&a &e_7^3d_4 &6 &25&de\cr
35&b &a_{13}d_{11} &2 &17&da\cr
\cr
36&a &a_{18}e_6 &2 &13&aa\cr
\cr
37&a*&a_{17}e_7a_1 &2 &17&aaA\cr
37&d*&d_{10}e_7^2a_1 &2 &25&ddD\cr
\cr
39&a &d_{12}e_7d_6 &1 &25&ddd\cr
\cr
45&d*&d_{12}^2a_1 &2 &25&dD\cr
\cr
47&a &d_{10}e_8e_7 &1 &25&edd\cr
47&b &d_{14}d_{10}a_1 &1 &25&ddd\cr
47&c &a_{17}e_8 &2 &17&da\cr
\cr
48&a &a_{23}a_2 &2 &13&ad\cr
\cr
51&a*&a_{24}a_1 &2 &13&aA\cr
\cr
61&d*&d_{16}e_8a_1 &1 &25&ddD\cr
61&e*&e_8^3a_1 &6 &25&eE\cr
\cr
63&a &d_{18}e_7 &1 &25&dd\cr
\cr
93&d*&d_{24}a_1 &1 &25&dD\cr
}

\proclaim Table $2$.~The norm $4$ vectors of $II_{1,25}$. 

There is a natural 1:1 correspondence between the elements of the
following sets:
\item{(1)} 
Orbits of norm $4$ vectors $u$ in $II_{1,25}$ under 
$\Aut(II_{1,25})$. 
\item{(2)}
Orbits of norm $4$ vectors in the fundamental domain $D$ of
$II_{1,25}$ under $\Aut(D)$.
\item{(3)}
Orbits of norm $1$ vectors $v$ of $I_{1,25}$ under $\Aut(I_{1,25})$. 
\item{(4)}
25 dimensional unimodular positive definite lattices $L$.
\item{(5)} 
Unimodular lattices $L_1$ of dimension at most 25 with no vectors of
norm 1.
\item{(6)}
25 dimensional even lattices $L_2$ of determinant 4. 

$L_1$ is the orthogonal complement of the norm 1 vectors of $L$, $L_2$
is the lattice of elements of $L$ of even norm, $L_2(-1)$ is isomorphic to
$u^\perp$, and $L(-1)$ is isomorphic to $v^\perp$. Table 2 lists the
665 elements of any of these sets.

The height  is the height of the norm $4$ vector $u$ of $D$, in
other words $(u,w)$ where $w$ is the Weyl vector of $D$. The things
in table 2 are listed in increasing order of their height.

Dim is the dimension of the lattice $L_1$. A capital $E$ after the
dimension means that $L_1$ is even.

The column ``roots'' gives the Dynkin diagram of the norm 2 vectors of
$L_2$ arranged into orbits under $\Aut(L_2)$.

``Group'' gives the order of the subgroup of $\Aut(D)$ fixing $u$. The
group $\Aut(L)\cong \Aut(L_2)$ is of the form $2\times R.G$ where $R$
is the group generated by the reflections of norm 2 vectors of $L$,
$G$ is the group described in the column ``group'', and 2 is the group
of order 2 generated by $-1$. If $\dim(L_1)\le 24$ then $\Aut(L_1)$ is
of the form $R.G$ where $R$ is the reflection group of $L_1$ and $G$
is as above.

For any root $r$ of $u^\perp$ the vector $v=u-r$ is a norm $2$ vector
of $II_{1,25}$. This vector $v$ can  be found as follows. Let $X$ be the
component of the Dynkin diagram of $u^\perp$ to which $u$ belongs and
let $h$ be the Coxeter number of $X$. Then $u-r$ is conjugate to a
norm $2$ vector of $II_{1,25}$ in $D$ of height $t-h+1$ (or $t-h$ if the
entry under ``Dim'' is $24E$) whose letter is the letter corresponding
to $X$ in the column headed ``norm $2$''. For example let $u$ be
the vector of height 6 and root system $a_2^2a_1^{10}$. Then the norm
$2$ vectors corresponding to roots from the components $a_2$ or $a_1$
have heights $6-3+1$ and $6-2+1$ and letters $a$ and $b$, so they are
the vectors $4a$ and $5b$ of table 1.

If $\dim(L_1)\le 24$ then the column ``neighbors'' gives the two even
neighbors of $L_1+I^{24-\dim(L_1)}$. If $\dim(L_1)\le 23$ then both
neighbors are isomorphic so only one is listed, and if $L_1$ is a
Niemeier lattice then the neighbor is preceded by 2 (to indicate that
the corresponding norm 0 vector is twice a primitive vector).  If the
two neighbors are isomorphic then there is an automorphism of $L$
exchanging them.

\bigskip
\halign{
\hfill$#$&~\hfill$#$&#&~\hfill$#$&~\hfill$#$&~\hfill#&~\hfill$#$&\hfill$#$\cr
\hbox{Height}&\hbox{Dim}&&\hbox{Roots}&\hbox{Group}&\hbox{norm $2$}&
\hbox{Neighbors}&\cr
1 &24&E &\hbox{None} &8315553613086720000 &&2\Lambda&      \cr
\cr
2 &23& &a_1^2 &84610842624000 &a&\Lambda &      \cr
2 &24& &\hbox{None} &1002795171840 &&\Lambda&A_1^{24}      \cr
\cr
3 &25& &a_1^2 &88704000 &a&      &      \cr
\cr
4 &24& &a_1^8 &20643840 &a&A_1^{24}    &A_1^{24}    \cr
4 &25& &a_2^2 &26127360 &a&      &      \cr
4 &25& &a_1^6 &138240 &a&      &      \cr
\cr
5 &24& &a_1^{12} &190080 &a& A_1^{24}   & A_2^{12}   \cr
5 &25& &a_2a_1^7 &5040 &aa&      &      \cr
5 &25& &a_1^{10} &1920 &a&      &      \cr
\cr
6 &23& &a_1^{16}a_1^2 &645120 &ca&  A_1^{24}  &      \cr
6 &24& &a_2^2a_1^{10} &5760 &ab&    A_2^{12}&A_2^{12}    \cr
6 &24& &a_1^{16} &43008 &c&    A_1^{24}& A_3^8   \cr
6 &25& &a_3a_1^8 &21504 &ac&      &      \cr
6 &25& &a_2^2a_1^8 &128 &ab&      &      \cr
6 &25& &a_2a_1^{10}a_1 &120 &abc&      &      \cr
6 &25& &a_1^8a_1^6 &1152 &bc&      &      \cr
\cr
7 &24& &a_2^4a_1^8 &384 &bb&    A_2^{12}&A_3^8    \cr
7 &24&E &a_1^{24} &244823040 &a&   2A_1^{24}&      \cr
7 &25& &a_2^5a_1^3 &720 &bb&      &      \cr
7 &25& &a_3a_2a_1^9 &72 &acb&      &      \cr
7 &25& &a_2^4a_1^4a_1^2 &24 &bba&      &      \cr
7 &25& &a_3a_1^{12} &1440 &ab&      &      \cr
7 &25& &a_2^3a_1^6a_1^3 &12 &bbb&      &      \cr
7 &25& &a_2^2a_1^{12} &144 &cb&      &      \cr
\cr
8 &22& &a_3a_1^{22} &887040 &ae&    A_1^{24}&      \cr
8 &23& &a_2^6a_1^6a_1^2 &1440 &bda&   A_2^{12} &      \cr
8 &24& &a_3^2a_1^{12} &768 &cc&    A_3^8& A_3^8   \cr
8 &24& &a_3a_2^4a_1^6 &96 &bbb&    A_3^8&A_3^8    \cr
8 &24& &a_2^8 &672 &a&    A_3^8&A_3^8    \cr
8 &24& &a_2^6a_1^6 &240 &bd&   A_2^{12} &A_4^6    \cr
8 &24& &a_1^{24} &138240 &e&    A_1^{24}&D_4^6   \cr
8 &25& &a_4a_1^{12} &1440 &ad&      &      \cr
8 &25& &a_3a_2^4a_1^4 &16 &bbb&      &      \cr
8 &25& &a_3^2a_1^8a_1^2 &64 &cbc&      &      \cr
8 &25& &a_3a_2^3a_1^6a_1 &12 &bbbc&      &      \cr
8 &25& &a_3a_2^3a_1^3a_1^3a_1 &6 &bbbbd&      &      \cr
8 &25& &a_2^4a_2^2a_1^4 &8 &bab&      &      \cr
8 &25& &a_3a_2^2a_1^4a_1^4a_1^2 &16 &bbcbd&      &      \cr
8 &25& &a_2^4a_2a_1^4a_1^2a_1 &8 &bbbdb&      &      \cr
8 &25& &a_2^4a_1^8a_1^2 &48 &bdc&      &      \cr
8 &25& &a_3a_1^{15}a_1 &720 &cce&      &      \cr
\cr
9 &24& &a_3^2a_2^4a_1^4 &16 &bbb&  A_3^8  &A_4^6    \cr
9 &24& &a_2^8a_1^4 &384 &dc&    A_2^{12}&A_5^4D_4    \cr
9 &25& &a_4a_2^3a_1^6a_1 &12 &bdbb&      &      \cr
9 &25& &a_3^2a_2^4a_1^2 &16 &bba&      &      \cr
9 &25& &a_3^2a_2^2a_2a_1^4a_1 &4 &bbcba&      &      \cr
9 &25& &a_3a_3a_2^2a_2a_1^2a_1^2a_1 &2 &bbbbbbb&      &      \cr
9 &25& &a_3a_2^6a_1^2 &6 &abb&      &      \cr
9 &25& &a_3^2a_2^2a_1^4a_1^4 &8 &bcbb&      &      \cr
9 &25& &a_3a_2^4a_2a_1^4a_1 &8 &bbbbc&      &      \cr
9 &25& &a_3a_2^2a_2^2a_2a_1^2a_1^2a_1 &2 &bbbdbbb&      &      \cr
\cr
10 &22& &a_3a_2^{10} &2880 &ac&   A_2^{12} &      \cr
10 &23& &a_3^4a_1^8a_1^2 &384 &cfa&  A_3^8  &      \cr
10 &23& &a_3^3a_2^4a_1^2a_1^2 &48 &bbea& A_3^8   &      \cr
10 &24& &a_3^4a_2^2a_1^2 &32 &bab&    A_4^6&A_4^6    \cr
10 &24& &a_4a_3a_2^4a_1^4 &16 &bdbc&    A_4^6& A_4^6   \cr
10 &24& &a_3^2a_3a_2^4a_1^2 &16 &bbbe&    A_3^8&A_5^4D_4    \cr
10 &24& &a_3^4a_1^4a_1^4 &48 &cdf&    A_3^8& A_5^4D_4   \cr
10 &24& &a_3^4a_1^8 &384 &cd&    A_3^8&D_4^6   \cr
10 &24&E &a_2^{12} &190080 &a&   2A_2^{12}&      \cr
10 &25& &a_3^5 &1920 &c&      &      \cr
10 &25& &d_4a_2^4a_1^6 &144 &bcd&      &      \cr
10 &25& &d_4a_3a_1^{12} &576 &ced&      &      \cr
10 &25& &a_5a_1^{15} &720 &cf&      &      \cr
10 &25& &a_4a_3a_2^4a_1^2 &8 &bdbb&      &      \cr
10 &25& &a_3^3a_3a_2a_1^3 &6 &bcbb&      &      \cr
10 &25& &a_4a_3a_2^2a_2a_1^2a_1^2a_1 &2 &bdbbbce&      &      \cr
10 &25& &a_3^3a_3a_1^4a_1^2 &24 &bcce&      &      \cr
10 &25& &a_3^2a_3^2a_1^4a_1^2 &16 &cbbd&      &      \cr
10 &25& &a_4a_2^6a_1^2 &6 &abe&      &      \cr
10 &25& &a_4a_3a_2^2a_1^4a_1^2a_1^2 &4 &bdbccf&      &      \cr
10 &25& &a_3^2a_3a_2^2a_2a_1^2a_1 &2 &bbbbbc&      &      \cr
10 &25& &a_3^2a_3a_2^2a_2a_1a_1a_1 &2 &bbbadbc&      &      \cr
10 &25& &a_4a_2^5a_1^5 &10 &bbc&      &      \cr
10 &25& &a_3^2a_2^6 &48 &da&      &      \cr
10 &25& &a_3^2a_2^4a_2^2 &8 &bba&      &      \cr
10 &25& &a_3^3a_2^2a_1^3a_1^2a_1 &12 &bbdcf&      &      \cr
10 &25& &a_3a_3a_3a_2^2a_1^2a_1a_1a_1a_1 &2 &cbbbcebfd&      &      \cr
10 &25& &a_3a_3a_2^2a_2^2a_2a_1^2a_1 &2 &bbbbbee&      &      \cr
10 &25& &a_3^2a_2^4a_1^4a_1^2 &8 &dbcf&      &      \cr
10 &25& &a_3^3a_1^{12} &48 &cf&      &      \cr
10 &25& &a_3a_2^6a_2a_1^3 &12 &dbce&      &      \cr
\cr
11 &24& &a_4a_3^2a_3a_2^2a_1^2 &4 &bbbbb&    A_4^6&A_5^4D_4    \cr
11 &24& &a_4^2a_2^4a_1^4 &16 &dca&    A_4^6& A_5^4D_4   \cr
11 &24& &a_3^6 &240 &a&    A_4^6&D_4^6   \cr
11 &24& &a_3^4a_2^4 &24 &be&   A_3^8 &A_6^4    \cr
11 &25& &a_5a_2^4a_2a_1^4 &8 &befb&      &      \cr
11 &25& &d_4a_3a_2^4a_1^4 &8 &bcda&      &      \cr
11 &25& &a_4a_3^2a_3a_2a_1^2a_1 &2 &bbbcba&      &      \cr
11 &25& &a_4^2a_2^2a_2a_1^4a_1 &4 &dcbbb&      &      \cr
11 &25& &a_4a_3^2a_2^4 &4 &bbb&      &      \cr
11 &25& &a_4a_3^3a_1^6 &6 &cbb&      &      \cr
11 &25& &a_4a_3^2a_2^2a_2a_1^2a_1 &2 &bbbfab&      &      \cr
11 &25& &a_4a_3a_3a_2a_2a_2a_1a_1a_1 &1 &bbbbccbba&      &      \cr
11 &25& &a_4a_3a_3a_2a_2a_2a_1a_1a_1 &1 &bbbbecbbb&      &      \cr
11 &25& &a_3^4a_3a_1^4 &8 &bab&      &      \cr
11 &25& &a_3^2a_3^2a_2^2a_2a_1 &2 &bbbbb&      &      \cr
11 &25& &a_3^2a_3a_3a_2^2a_2a_1 &2 &bbabda&      &      \cr
11 &25& &a_4a_3a_2^2a_2^2a_2a_1^2a_1 &2 &bbeccba&      &      \cr
11 &25& &a_3^2a_3^2a_2^2a_1^4 &4 &bbdb&      &      \cr
\cr
12 &22& &a_3^6a_3a_1^2 &96 &dag&   A_3^8 &      \cr
12 &23& &a_4^2a_3^2a_2^2a_1^2a_1^2 &8 &bcbha& A_4^6   &      \cr
12 &23& &a_4a_3^5a_1^2 &40 &aba&    A_4^6&      \cr
12 &24& &d_4a_4a_2^6 &24 &dca&    A_5^4D_4& A_5^4D_4   \cr
12 &24& &d_4a_3^4a_1^4 &32 &cde&    A_5^4D_4& A_5^4D_4   \cr
12 &24& &a_5a_3^3a_1^6a_1 &24 &cfec&    A_5^4D_4&A_5^4D_4    \cr
12 &24& &a_4^2a_3^2a_3a_1^2 &8 &bbbd&    A_5^4D_4& A_5^4D_4   \cr
12 &24& &a_5a_3^2a_2^4a_1 &16 &bebf&    A_5^4D_4& A_5^4D_4   \cr
12 &24& &d_4a_3^4a_1^4 &48 &cdc&    D_4^6 &A_5^4D_4  \cr
12 &24& &d_4^2a_1^{16} &1152 &eb&    D_4^6&D_4^6  \cr
12 &24& &a_4^2a_3^2a_2^2a_1^2 &4 &bcbh&   A_4^6 &A_6^4    \cr
12 &24& &a_3^4a_3^2a_1^4 &32 &fdg&    A_3^8&A_7^2D_5^2    \cr
12 &25& &a_5a_3^2a_3a_1^4a_1 &8 &cefdc&      &      \cr
12 &25& &a_5a_3^2a_2^2a_2a_1a_1 &2 &bebbdh&      &      \cr
12 &25& &d_4a_3^2a_3a_2^2a_1^2 &4 &bddac&      &      \cr
12 &25& &d_4a_4a_2^4a_1^4 &8 &dcae&      &      \cr
12 &25& &a_5a_3^2a_2^2a_1^2a_1^2a_1 &4 &bfbdhf&      &      \cr
12 &25& &a_5a_3a_3a_2^2a_1^2a_1a_1a_1 &2 &bfebhdee&      &      \cr
12 &25& &a_4^2a_3a_3a_2a_1^2a_1 &2 &bbcade&      &      \cr
12 &25& &a_4a_4a_3a_3a_2a_1a_1a_1 &1 &bbccbddd&      &      \cr
12 &25& &a_4^2a_3a_2^4 &4 &bbb&      &      \cr
12 &25& &a_4^2a_3^2a_1^4a_1^2 &4 &bche&      &      \cr
12 &25& &d_4a_3^2a_2^4a_1a_1 &8 &bdaeg&      &      \cr
12 &25& &d_4a_3^3a_1^6a_1a_1 &12 &cdegb&      &      \cr
12 &25& &a_4a_3^2a_3^2a_2a_1 &2 &abcbc&      &      \cr
12 &25& &a_4^2a_3a_2^2a_2a_1a_1a_1 &2 &bcbbfdh&      &      \cr
12 &25& &a_4a_4a_3a_2a_2a_2a_1a_1a_1 &1 &bbcabbhhd&      &      \cr
12 &25& &a_4a_3^4a_1^4 &8 &ach&      &      \cr
12 &25& &a_4a_3^2a_3a_3a_1^2a_1^2 &2 &bbcfde&      &      \cr
12 &25& &a_4a_3^3a_2^3a_1 &6 &bbag&      &      \cr
12 &25& &a_4a_3a_3a_3a_2a_2a_2a_1 &1 &bbbebbah&      &      \cr
12 &25& &a_3^4a_3a_3a_1^2 &8 &bfdc&      &      \cr
12 &25& &a_4a_3^2a_3a_2^2a_1^2a_1a_1 &2 &bccbhge&      &      \cr
12 &25& &a_4a_3^2a_3a_2^2a_1^2a_1^2 &2 &bcfbhh&      &      \cr
12 &25& &a_3^2a_3^2a_3a_2^2a_1^2 &4 &edbbh&      &      \cr
\cr
13 &24& &a_4^4a_1^4 &24 &cc&    A_5^4D_4&A_6^4    \cr
13 &24& &a_5a_4a_3a_3a_2a_2a_1 &2 &bebbddd&  A_5^4D_4  & A_6^4   \cr
13 &24& &a_4^2a_3^4 &16 &bb&   A_4^6 &A_7^2D_5^2    \cr
13 &24& &a_4^2a_4a_3a_2^2a_1^2 &4 &ccahb&  A_4^6  & A_7^2D_5^2   \cr
13 &24&E &a_3^8 &2688 &a& 2A_3^8  &      \cr
13 &25& &a_6a_2^6a_1^3 &12 &dhd&      &      \cr
13 &25& &a_4^4a_1^2 &24 &ca&      &      \cr
13 &25& &a_5a_4a_3^2a_2a_1^2 &2 &bfbdc&      &      \cr
13 &25& &d_4a_4a_3^3a_1^2 &6 &bdac&      &      \cr
13 &25& &d_4^2a_2^6 &72 &cc&      &      \cr
13 &25& &a_5a_4a_3a_2a_2a_2a_1a_1 &1 &bebhdeda&      &      \cr
13 &25& &a_5a_4a_3a_2a_2a_2a_1a_1 &1 &bebhdddd&      &      \cr
13 &25& &d_4a_4a_3a_3a_2a_2a_1a_1 &1 &bdaaeccb&      &      \cr
13 &25& &a_4^2a_4a_3a_2^2 &2 &bcbd&      &      \cr
13 &25& &a_4a_4a_4a_3a_2a_1a_1a_1 &1 &bccbdbcd&      &      \cr
13 &25& &a_5a_3^3a_2^3 &6 &bbd&      &      \cr
13 &25& &a_5a_3^2a_3a_2^2a_2 &2 &abbhc&      &      \cr
13 &25& &a_5a_4a_2^2a_2^2a_2a_1^2 &2 &behdfd&      &      \cr
13 &25& &a_5a_3^2a_3a_2^2a_1^2a_1 &2 &bbbedd&      &      \cr
13 &25& &a_4a_4a_3a_3a_3a_2a_1 &1 &bcbbbdc&      &      \cr
13 &25& &a_4a_4a_3a_3a_3a_2a_1 &1 &bbbabdb&      &      \cr
13 &25& &a_4^3a_2^3a_1^3 &6 &cdb&      &      \cr
13 &25& &d_4a_3^3a_2^3a_2 &6 &bacg&      &      \cr
13 &25& &a_4^2a_3a_3a_2^2a_2a_1 &2 &ebaddd&      &      \cr
13 &25& &a_5a_2^9 &72 &cf&      &      \cr
\cr
14 &21& &d_4a_3^7 &336 &ag& A_3^8   &      \cr
14 &22& &a_4^4a_3a_2^2 &16 &aab& A_4^6   &      \cr
14 &23& &a_5a_4^2a_3a_3a_1^2a_1 &4 &bbddaf& A_5^4D_4   &      \cr
14 &23& &d_4a_4^3a_2^2a_1^2 &12 &caca&   A_5^4D_4 &      \cr
14 &23& &a_5d_4a_3^3a_1^3a_1^2 &12 &dfega& A_5^4D_4   &      \cr
14 &23& &a_5^2a_3a_2^4a_1^2 &16 &efda&   A_5^4D_4 &      \cr
14 &23& &d_4^2a_3^4a_1^2 &96 &dcd&   D_4^6&      \cr
14 &24& &a_5a_4^3a_1^3 &12 &cbe&   A_6^4 & A_6^4   \cr
14 &24& &a_5^2a_3^2a_2^2 &8 &eda&    A_6^4&A_6^4    \cr
14 &24& &a_6a_3^3a_2^3 &12 &bhd&    A_6^4& A_6^4   \cr
14 &24& &d_4a_4^2a_4a_2^2 &4 &caac&   A_5^4D_4 & A_7^2D_5^2   \cr
14 &24& &a_5d_4a_3^2a_3a_1^2a_1 &4 &dfecbg&  A_5^4D_4  & A_7^2D_5^2   \cr
14 &24& &a_5^2a_3^2a_1^2a_1^2a_1^2 &8 &febcg&    A_5^4D_4& A_7^2D_5^2   \cr
14 &24& &a_5a_4^2a_3a_3a_1 &4 &bbddf&    A_5^4D_4&A_7^2D_5^2    \cr
14 &24& &d_4^2a_3^4 &32 &dc&   D_4^6& A_7^2D_5^2   \cr
14 &24& &a_4^3a_3^3 &12 &bh&    A_4^6&A_8^3    \cr
14 &24& &a_3^8 &384 &g&    A_3^8& D_6^4  \cr
14 &25& &d_5a_3^2a_2^4a_1^2 &8 &bgbb&      &      \cr
14 &25& &d_5a_3^3a_1^8 &48 &cgc&      &      \cr
14 &25& &a_6a_3^2a_3a_2^2a_1 &2 &bhhcf&      &      \cr
14 &25& &a_5^2a_3a_3a_1^4 &8 &efee&      &      \cr
14 &25& &a_6a_3^2a_3a_2a_1^2a_1^2 &2 &bhhdeg&      &      \cr
14 &25& &d_4a_4^3a_1^3a_1 &6 &cabc&      &      \cr
14 &25& &a_5d_4a_3a_3a_2^2a_1 &2 &deeccb&      &      \cr
14 &25& &a_5a_5a_3a_2^2a_1^2a_1a_1 &2 &efddefg&      &      \cr
14 &25& &a_5a_4^2a_3a_2a_1^2 &2 &bbfdf&      &      \cr
14 &25& &a_5a_4a_4a_3a_2a_1a_1 &1 &cbbdabe&      &      \cr
14 &25& &d_4a_4^2a_3^2a_1^2 &4 &baeb&      &      \cr
14 &25& &d_4^2a_3^2a_3a_1^4 &8 &dcbb&      &      \cr
14 &25& &a_5a_4^2a_3a_1^2a_1^2a_1 &2 &cbhcee&      &      \cr
14 &25& &a_5d_4a_3^2a_1^4a_1^2a_1 &4 &dfegcg&      &      \cr
14 &25& &a_5a_4a_3a_3a_3a_2 &1 &bbhddc&      &      \cr
14 &25& &a_5a_4a_3a_3a_3a_2 &1 &bbeddc&      &      \cr
14 &25& &a_5a_4^2a_2^2a_2a_1^2 &2 &cbddf&      &      \cr
14 &25& &a_4^4a_2^2 &8 &ba&      &      \cr
14 &25& &d_4a_4a_4a_3a_2a_2a_1a_1 &1 &caaecbbg&      &      \cr
14 &25& &a_5a_4a_3a_3a_3a_1a_1a_1 &1 &bbhedbfg&      &      \cr
14 &25& &a_5a_3^2a_3^2a_3a_1 &4 &dddcg&      &      \cr
14 &25& &a_5a_3^2a_3a_3a_3a_1 &2 &bheddf&      &      \cr
14 &25& &a_5a_4a_3^2a_2^2a_1^2a_1 &2 &cbhdge&      &      \cr
14 &25& &a_4^2a_4a_3^2a_2a_1 &2 &bbdbf&      &      \cr
14 &25& &a_4a_4a_4a_3a_3a_2a_1 &1 &abbdhcf&      &      \cr
14 &25& &a_4^2a_4a_3a_2^2a_2a_1 &2 &bahbdf&      &      \cr
14 &25& &d_4a_3^4a_3a_1^4 &8 &fegg&      &      \cr
\cr
15 &24& &a_5^2a_4^2a_1^2 &4 &bdc&  A_6^4  & A_7^2D_5^2   \cr
15 &24& &a_6a_4a_4a_3a_2a_1a_1 &2 &chhceaa& A_6^4   & A_7^2D_5^2   \cr
15 &24& &a_5^2a_4a_3a_2^2 &4 &bfdf&   A_5^4D_4 & A_8^3   \cr
15 &24& &a_4^4a_3^2 &16 &hb&    A_4^6& A_9^2D_6   \cr
15 &25& &d_5a_3^5 &20 &ab&      &      \cr
15 &25& &d_5a_4a_3^2a_2^2a_1^2 &2 &bgbba&      &      \cr
15 &25& &a_6a_4^2a_3a_1^2a_1 &2 &chdcb&      &      \cr
15 &25& &a_6a_4a_4a_2a_2a_1a_1a_1 &1 &chhefacc&      &      \cr
15 &25& &a_5d_4a_4^2a_1^2a_1 &2 &abeab&      &      \cr
15 &25& &a_5a_5a_4a_3a_2a_1 &1 &bbddea&      &      \cr
15 &25& &a_6a_4a_3a_3a_2a_2a_1 &1 &bhdcfga&      &      \cr
15 &25& &a_6a_4a_3^2a_2^2a_1 &2 &bhdfc&      &      \cr
15 &25& &d_4^2a_4^2a_2^2 &4 &acb&      &      \cr
15 &25& &a_5^2a_4a_3a_1^2a_1^2 &2 &bedcc&      &      \cr
15 &25& &a_5d_4a_4a_3a_2a_2a_1 &1 &abecbba&      &      \cr
15 &25& &a_5^2a_3^2a_3a_1^2 &2 &bdac&      &      \cr
15 &25& &a_5a_4^2a_4a_2a_1^2 &2 &adfec&      &      \cr
15 &25& &a_5a_4^2a_4a_2a_1^2 &2 &addca&      &      \cr
15 &25& &a_5a_4a_4a_4a_2a_1a_1 &1 &bhddeca&      &      \cr
15 &25& &a_4^5 &5 &d&      &      \cr
15 &25& &a_5a_4a_4a_3a_3a_2 &1 &bhdcab&      &      \cr
15 &25& &d_4a_4^2a_3^2a_3 &2 &bccb&      &      \cr
\cr
16 &21& &d_4a_4^5 &40 &ab&   A_4^6 &      \cr
16 &22& &a_5^2d_4a_3^2a_3 &8 &ceba&    A_5^4D_4&      \cr
16 &22& &a_5^3a_3a_3a_1^3 &12 &ecad&   A_5^4D_4 &      \cr
16 &22& &d_4^4a_3a_1^6 &144 &bdc&   D_4^6&      \cr
16 &23& &a_6a_5a_4a_3a_2a_1^2a_1 &2 &bhdecah&  A_6^4  &      \cr
16 &23& &a_5^3a_4a_1^2a_1 &6 &daag&    A_6^4&      \cr
16 &24& &d_5d_4a_3^4 &16 &dgb&  A_7^2D_5^2  &   A_7^2D_5^2 \cr
16 &24& &d_5a_5a_3^2a_3a_1^2a_1 &4 &fgbceb&  A_7^2D_5^2  &  A_7^2D_5^2  \cr
16 &24& &a_5^2d_4^2a_1^2 &16 &cef&  A_7^2D_5^2  & A_7^2D_5^2   \cr
16 &24& &a_7a_3^4a_1^4 &16 &fge&  A_7^2D_5^2  &  A_7^2D_5^2  \cr
16 &24& &d_5a_4a_4^2a_2^2 &4 &cbba&  A_7^2D_5^2  & A_7^2D_5^2   \cr
16 &24& &a_6d_4a_4^2a_2 &4 &ahcb&  A_7^2D_5^2  &  A_7^2D_5^2  \cr
16 &24& &a_6a_5a_4a_3a_2a_1 &2 &bhdech&  A_6^4  &  A_8^3  \cr
16 &24& &a_5^2d_4a_3^2a_1^2 &4 &eegd&   A_5^4D_4 & A_9^2D_6  \cr
16 &24& &a_5a_5a_4^2a_3 &4 &dddf&    A_5^4D_4&  A_9^2D_6  \cr
16 &24& &a_5^2d_4a_3^2a_1^2 &8 &eebd&   A_5^4D_4 &  D_6^4 \cr
16 &24& &d_4^4a_1^8 &48 &bc&  D_4^6 & D_6^4  \cr
16 &24&E &a_4^6 &240 &a&  2A_4^6 &      \cr
16 &25& &a_7a_3^4a_1^2 &8 &fff&      &      \cr
16 &25& &a_7a_3a_3a_3a_2^2a_1^2 &2 &efggch&      &      \cr
16 &25& &a_5d_4^3a_1^3 &36 &bcb&      &      \cr
16 &25& &d_5a_4^2a_3a_3a_1a_1 &2 &bbcbdb&      &      \cr
16 &25& &d_5d_4a_3^3a_1^3a_1 &6 &dgbec&      &      \cr
16 &25& &d_5a_5a_3a_2^4a_1 &4 &egcad&      &      \cr
16 &25& &a_6a_5a_4a_2a_2a_1a_1 &1 &bhdcagh&      &      \cr
16 &25& &d_5a_4a_3^4 &4 &bbb&      &      \cr
16 &25& &d_5a_4^2a_3a_2^2a_1^2 &2 &cbcae&      &      \cr
16 &25& &a_6d_4a_4a_3a_2a_1a_1 &1 &ahcgafe&      &      \cr
16 &25& &a_5^2a_5a_3a_2 &2 &deeb&      &      \cr
16 &25& &a_6a_5a_3a_3a_2a_2 &1 &bhfebc&      &      \cr
16 &25& &a_6a_5a_3a_3a_2a_1a_1a_1 &1 &bhegchhh&      &      \cr
16 &25& &a_6a_4^2a_3^2a_1 &2 &bcee&      &      \cr
16 &25& &a_5a_5a_5a_2^2a_1^2a_1 &2 &defchd&      &      \cr
16 &25& &a_5^2d_4a_3a_2^2 &4 &cfcb&      &      \cr
16 &25& &a_5a_5d_4a_3a_2^2 &2 &ccdga&      &      \cr
16 &25& &a_6d_4a_3^2a_2^2a_2 &2 &ahgab&      &      \cr
16 &25& &a_5a_5a_4a_4a_2a_1 &1 &ddddcg&      &      \cr
16 &25& &a_5a_5a_4a_4a_2a_1 &1 &ddadad&      &      \cr
16 &25& &a_6a_4a_4a_3a_2a_2a_1 &1 &bddecah&      &      \cr
16 &25& &a_5d_4a_4^2a_3a_1 &2 &edccb&      &      \cr
16 &25& &a_5d_4^2a_3^2a_1a_1a_1 &2 &cebecf&      &      \cr
16 &25& &a_5^2a_4a_3^2a_1^2 &2 &haeh&      &      \cr
16 &25& &a_5^2a_4a_3a_3a_1^2 &2 &ddebd&      &      \cr
16 &25& &a_5a_5a_4a_3a_3a_1a_1 &1 &eddefdg&      &      \cr
16 &25& &a_5^2a_3^4 &8 &cf&      &      \cr
16 &25& &a_5^2a_4a_3a_2^2a_1^2 &2 &hdcch&      &      \cr
16 &25& &d_4a_4^2a_4a_3^2 &2 &hcbb&      &      \cr
\cr
17 &24& &a_7a_4^2a_3^2 &4 &bfc&   A_7^2D_5^2 &  A_8^3  \cr
17 &24& &a_6^2a_4a_3a_1^2 &4 &hebb& A_7^2D_5^2   &A_8^3    \cr
17 &24& &a_6a_5a_5a_3a_2 &2 &dddcg&   A_6^4 &  A_9^2D_6  \cr
17 &24& &a_6^2a_3^2a_2^2 &4 &hah&  A_6^4  &  A_9^2D_6  \cr
17 &24& &a_5^4 &24 &a& A_6^4   & D_6^4  \cr
17 &25& &a_7a_4a_4a_3a_2a_1 &1 &bfgcgb&      &      \cr
17 &25& &a_6^2a_3^2a_2a_1 &2 &hcfa&      &      \cr
17 &25& &a_7a_4^2a_2^2a_2a_1 &2 &bfhea&      &      \cr
17 &25& &d_5d_4a_4^2a_2^2 &2 &abbb&      &      \cr
17 &25& &d_5a_5a_4a_3a_3a_1 &1 &bbbaab&      &      \cr
17 &25& &d_5a_5a_4a_3a_2a_2a_1 &1 &bbbaedb&      &      \cr
17 &25& &a_6a_5d_4a_3a_2a_1 &1 &ecdafb&      &      \cr
17 &25& &a_6a_5a_4^2a_1^2 &2 &ddeb&      &      \cr
17 &25& &a_6a_5a_4a_3a_3 &1 &dcgbc&      &      \cr
17 &25& &a_6a_5a_4a_3a_3 &1 &dcfcc&      &      \cr
17 &25& &a_6a_5a_4a_3a_3 &1 &dceac&      &      \cr
17 &25& &a_5^2d_4a_4a_3 &2 &cdbb&      &      \cr
17 &25& &a_6a_4a_4a_4a_3a_1 &1 &deffab&      &      \cr
\cr
18 &21& &a_5^3d_4d_4a_1 &12 &bcab&  A_5^4D_4  &      \cr
18 &22& &a_6^2a_4^2a_3 &4 &caa&  A_6^4  &      \cr
18 &23& &a_6d_5a_4a_4a_2a_1^2 &2 &bhaaba&  A_7^2D_5^2  &      \cr
18 &23& &a_6^2d_4a_4a_1^2 &4 &ceba&  A_7^2D_5^2  &      \cr
18 &23& &d_5a_5^2d_4a_1^2a_1^2 &4 &ebcfa&  A_7^2D_5^2  &      \cr
18 &23& &a_7a_5a_4^2a_1^2a_1 &4 &dfcag&  A_7^2D_5^2  &      \cr
18 &23& &a_7d_4^2a_3^2a_1^2 &8 &cgfa&  A_7^2D_5^2  &      \cr
18 &23& &d_5^2a_3^4a_1^2 &16 &gea&   A_7^2D_5^2 &      \cr
18 &24& &a_7a_5^2a_2^2 &8 &ffa&   A_8^3 & A_8^3   \cr
18 &24& &a_7a_5a_4^2a_1 &4 &dfcg&   A_7^2D_5^2 & A_9^2D_6   \cr
18 &24& &a_7a_5d_4a_3a_1a_1a_1 &2 &eggfdfc&A_7^2D_5^2    &A_9^2D_6    \cr
18 &24& &a_6d_5a_4a_4a_2 &2 &bhaab&  A_7^2D_5^2  &  A_9^2D_6  \cr
18 &24& &d_5^2a_3^4 &16 &gb&  A_7^2D_5^2  & D_6^4  \cr
18 &24& &d_5a_5^2d_4a_1^2 &4 &ebbc& A_7^2D_5^2   &  D_6^4 \cr
18 &24& &a_5^2a_5d_4a_3a_1 &4 &gbcdb& A_5^4D_4   &  A_{11}D_7E_6 \cr
18 &24& &a_5^4a_1^4 &48 &cb&   A_5^4D_4 &E_6^4  \cr
18 &25& &d_6a_3^4a_3a_1^2 &8 &ddcd&      &      \cr
18 &25& &a_7a_5^2a_1^4 &8 &fge&      &      \cr
18 &25& &a_7a_5d_4a_2^2a_1 &2 &egfbc&      &      \cr
18 &25& &a_7d_4^2a_3a_1^4 &4 &cgef&      &      \cr
18 &25& &a_7a_5a_4a_3a_2 &1 &dfchb&      &      \cr
18 &25& &a_6^2a_5a_2a_1a_1 &2 &deaed&      &      \cr
18 &25& &a_7a_5a_4a_3a_1a_1a_1 &1 &dgchfeg&      &      \cr
18 &25& &a_6d_5a_4a_3a_2a_1a_1 &1 &bhaebfd&      &      \cr
18 &25& &d_5a_5a_5a_4a_1a_1 &1 &dbcacc&      &      \cr
18 &25& &a_7a_5a_3^2a_3a_1 &2 &dggfd&      &      \cr
18 &25& &a_7a_5a_3^2a_3a_1 &2 &dfhfg&      &      \cr
18 &25& &a_6a_6a_4a_3a_2a_1 &1 &cdchbg&      &      \cr
18 &25& &a_6a_5a_5a_4a_1 &1 &aeecc&      &      \cr
18 &25& &d_5a_5d_4a_3^2a_1^2a_1 &2 &ebcefb&      &      \cr
18 &25& &d_5d_4^2a_3^2a_3 &4 &cbbc&      &      \cr
18 &25& &d_5a_5a_4^2a_3a_1 &2 &dbadf&      &      \cr
18 &25& &d_5a_5a_4^2a_3a_1 &2 &dbabb&      &      \cr
18 &25& &a_6a_5d_4a_4a_2a_1 &1 &cgeabf&      &      \cr
18 &25& &a_6a_5a_4a_4a_3 &1 &cfacg&      &      \cr
18 &25& &a_5^2d_4^2a_3a_1^2 &2 &bgff&      &      \cr
18 &25& &a_6a_4^2a_4^2a_1 &2 &bcag&      &      \cr
\cr
19 &24& &a_8a_4^2a_3^2 &4 &hhb&   A_8^3 &   A_9^2D_6 \cr
19 &24& &a_7a_6a_5a_2a_1 &2 &dfceb& A_8^3   &   A_9^2D_6 \cr
19 &24& &a_6a_6a_5a_4a_1 &2 &eeaha&  A_6^4  &A_{11}D_7E_6   \cr
19 &24&E &d_4^6 &2160 &d&  2D_4^6&      \cr
19 &24&E &a_5^4d_4 &48 &aa& 2A_5^4D_4  &      \cr
19 &25& &d_6a_4^2a_4a_2^2 &2 &addc&      &      \cr
19 &25& &a_8a_4a_4a_3a_2a_1 &1 &hghbfb&      &      \cr
19 &25& &a_8a_4^2a_2^2a_2a_1 &2 &hhgeb&      &      \cr
19 &25& &a_7a_6a_4a_2a_2a_1 &1 &dfheeb&      &      \cr
19 &25& &d_5^2a_4a_4a_2^2 &2 &bbec&      &      \cr
19 &25& &a_6d_5d_4a_4a_2 &1 &bcaec&      &      \cr
19 &25& &a_6d_5a_5a_3a_2a_1 &1 &bdabca&      &      \cr
19 &25& &a_7a_5^2a_3a_1^2 &2 &dcab&      &      \cr
19 &25& &d_5a_5^3a_1 &3 &aaa&      &      \cr
19 &25& &a_7a_5^2a_2^2a_1^2 &2 &dcdb&      &      \cr
19 &25& &a_7a_5a_4a_4a_2 &1 &achgc&      &      \cr
19 &25& &a_7d_4a_4a_4a_3 &1 &ccefb&      &      \cr
19 &25& &a_6^2a_5a_3a_2 &2 &fbad&      &      \cr
19 &25& &a_6a_5^3 &6 &cb&      &      \cr
19 &25& &a_6a_5^2a_5 &2 &ecb&      &      \cr
19 &25& &a_6a_5a_5d_4a_2 &1 &babcc&      &      \cr
19 &25& &d_5a_5^2a_4a_2^2 &2 &dadc&      &      \cr
\cr
20 &20& &d_5a_5^4 &16 &ad&  A_5^4D_4  &      \cr
20 &20& &d_5d_4^5 &120 &dc& D_4^6  &      \cr
20 &21& &a_6^3d_4a_2 &6 &aaa&  A_6^4  &      \cr
20 &22& &a_7d_5a_5a_3a_3a_1 &2 &bgecad& A_7^2D_5^2   &      \cr
20 &22& &d_5^2a_5^2a_3 &8 &bba&  A_7^2D_5^2  &      \cr
20 &22& &a_7^2a_3^2a_3a_1^2 &8 &gdae& A_7^2D_5^2   &      \cr
20 &23& &a_7a_6^2a_1^2a_1^2 &4 &ecga& A_8^3   &      \cr
20 &23& &a_7^2a_4a_3a_1^2 &4 &faea&  A_8^3  &      \cr
20 &23& &a_8a_5^2a_2^2a_1^2 &4 &dhba&  A_8^3  &      \cr
20 &24& &a_7^2d_4a_1^4 &8 &gff&  A_9^2D_6  &   A_9^2D_6 \cr
20 &24& &d_6a_5^2a_3^2 &8 &edd&  A_9^2D_6  &   A_9^2D_6 \cr
20 &24& &a_8a_5^2a_3 &4 &dge&  A_9^2D_6  & A_9^2D_6   \cr
20 &24& &d_6a_5^2a_3^2 &4 &edc&  D_6^4 &A_9^2D_6    \cr
20 &24& &d_6d_4^3a_1^6 &12 &bcb&  D_6^4 &D_6^4   \cr
20 &24& &a_7d_5a_5a_3a_1a_1a_1 &2 &cgefecd& A_7^2D_5^2   &A_{11}D_7E_6   \cr
20 &24& &d_5d_5a_5^2a_1^2 &4 &bced& A_7^2D_5^2   &  A_{11}D_7E_6 \cr
20 &24& &a_7d_5d_4a_3a_3 &2 &bgecf& A_7^2D_5^2   &  A_{11}D_7E_6 \cr
20 &24& &a_6a_6d_5a_4 &2 &aaeb&   A_7^2D_5^2 & A_{11}D_7E_6  \cr
20 &24& &d_5^2a_5^2a_1^2 &8 &cbc&  A_7^2D_5^2  &E_6^4  \cr
20 &24& &a_6^2a_5^2 &4 &ch&  A_6^4  &  A_{12}^2 \cr
20 &24& &d_4^6 &48 &c&  D_4^6 & D_8^3  \cr
20 &25& &a_8a_5^2a_1^2a_1^2 &2 &dhfg&      &      \cr
20 &25& &d_6a_5^2a_3a_1^2a_1a_1 &2 &eddfeb&      &      \cr
20 &25& &d_6a_5d_4a_3^2a_1 &2 &ccdcd&      &      \cr
20 &25& &a_7a_7a_3a_2^2a_1^2 &2 &fggbg&      &      \cr
20 &25& &a_8d_4a_4a_3a_3 &1 &chbff&      &      \cr
20 &25& &d_5^2a_5d_4a_1^2a_1 &2 &bbecb&      &      \cr
20 &25& &a_7d_5a_4^2a_1a_1 &2 &cfbcf&      &      \cr
20 &25& &a_7a_6a_5a_3 &1 &ecfe&      &      \cr
20 &25& &a_6^2d_5a_3a_1^2 &2 &aedd&      &      \cr
20 &25& &a_7d_5a_4a_3a_3 &1 &bfbfc&      &      \cr
20 &25& &a_7a_6a_5a_2a_1a_1a_1 &1 &echbgfe&      &      \cr
20 &25& &a_7a_6a_4a_4a_1 &1 &ecbad&      &      \cr
20 &25& &a_6^2a_6a_3a_1 &2 &abee&      &      \cr
20 &25& &a_6a_6a_6a_3a_1 &1 &caceg&      &      \cr
20 &25& &a_6d_5a_5a_4a_2a_1 &1 &aeebad&      &      \cr
20 &25& &a_6a_6a_5d_4a_1 &1 &abfhd&      &      \cr
20 &25& &a_7a_5^2a_4a_1^2 &2 &ehag&      &      \cr
20 &25& &d_5a_5a_5a_4^2 &2 &fbdb&      &      \cr
\cr
21 &24& &a_8a_6a_5a_2a_1 &2 &ehbga&  A_8^3  & A_{11}D_7E_6  \cr
21 &24& &a_7^2a_4^2 &4 &cg& A_7^2D_5^2   & A_{12}^2  \cr
21 &25& &d_6a_6a_4a_4a_2 &1 &cdcdd&      &      \cr
21 &25& &a_8a_6a_4a_3a_1 &1 &eggbb&      &      \cr
21 &25& &a_7a_6d_5a_2a_1a_1 &1 &aecfba&      &      \cr
21 &25& &a_6d_5d_5a_4a_2 &1 &baacc&      &      \cr
21 &25& &a_7d_5a_5a_4a_1 &1 &acbdb&      &      \cr
21 &25& &a_7a_6a_5a_4 &1 &cgbe&      &      \cr
\cr
22 &21& &a_7d_5^2d_4a_3 &4 &beac& A_7^2D_5^2   &      \cr
22 &22& &a_8a_6^2a_3 &4 &bba&   A_8^3 &      \cr
22 &23& &a_8a_7a_5a_1^2a_1 &2 &cgeac&  A_9^2D_6  &      \cr
22 &23& &a_8a_6d_5a_2a_1^2 &2 &abhba&   A_9^2D_6 &      \cr
22 &23& &a_7d_6a_5a_3a_1^2a_1 &2 &dgdfab& A_9^2D_6   &      \cr
22 &23& &a_9a_5a_4^2a_1^2 &4 &fgba& A_9^2D_6   &      \cr
22 &23& &d_6d_5a_5^2a_1^2 &4 &bdcd&  D_6^4 &      \cr
22 &23& &d_5^4a_1^2 &48 &bd&  D_6^4 &      \cr
22 &24& &a_9a_5d_4a_3a_1a_1 &2 &gfffeb& A_9^2D_6   & A_{11}D_7E_6  \cr
22 &24& &a_7d_6a_5a_3a_1 &2 &dgcfe&  A_9^2D_6  & A_{11}D_7E_6  \cr
22 &24& &a_8a_6d_5a_2 &2 &abhb& A_9^2D_6   &A_{11}D_7E_6   \cr
22 &24& &d_6d_5a_5^2 &4 &bdc&  D_6^4 & A_{11}D_7E_6  \cr
22 &24& &d_5^4 &48 &b& D_6^4  & E_6^4 \cr
22 &24& &a_8a_7a_4a_3 &2 &chbg&  A_8^3  & A_{12}^2  \cr
22 &24& &a_7d_5^2a_3^2 &4 &eed& A_7^2D_5^2   & D_8^3  \cr
22 &24& &a_7^2d_4^2 &8 &fd&  A_7^2D_5^2  & D_8^3  \cr
22 &24&E &a_6^4 &24 &a& 2A_6^4  &      \cr
22 &25& &e_6d_4^2a_3^3 &12 &cbc&      &      \cr
22 &25& &e_6a_5a_4^2a_3a_1 &2 &dbace&      &      \cr
22 &25& &e_6a_5^2a_2^4 &8 &fba&      &      \cr
22 &25& &d_7a_3^6 &24 &ge&      &      \cr
22 &25& &a_9a_5d_4a_2^2 &2 &gfgb&      &      \cr
22 &25& &a_9a_5a_4a_3a_1a_1 &1 &ffbgbc&      &      \cr
22 &25& &a_9a_5a_3^2a_3 &2 &fggf&      &      \cr
22 &25& &a_8a_7a_4a_2a_1 &1 &chbbc&      &      \cr
22 &25& &d_6d_5a_5a_3a_3a_1 &1 &bdcdbe&      &      \cr
22 &25& &a_7d_6a_3^2a_3a_1^2 &2 &dgfeb&      &      \cr
22 &25& &a_8d_5a_5a_3a_1 &1 &agffe&      &      \cr
22 &25& &a_8d_5a_5a_2a_2a_1 &1 &ahfbab&      &      \cr
22 &25& &a_8a_6a_5a_3 &1 &bbgg&      &      \cr
22 &25& &a_8a_6a_5a_3 &1 &cbef&      &      \cr
22 &25& &d_6a_5^2d_4a_3 &2 &bcdb&      &      \cr
22 &25& &d_6a_5^3a_1^3 &6 &cdb&      &      \cr
22 &25& &a_8a_5^2d_4 &2 &bfe&      &      \cr
22 &25& &a_8a_6a_4a_4a_1 &1 &cbbbc&      &      \cr
22 &25& &a_7a_7a_5a_3a_1 &1 &ghefc&      &      \cr
22 &25& &d_5^3a_3^3 &6 &ec&      &      \cr
22 &25& &a_7d_5d_4^2a_3 &2 &bffc&      &      \cr
\cr
23 &24& &a_8^2a_3^2 &4 &hb&  A_9^2D_6  &  A_{12}^2 \cr
23 &24& &a_9a_6a_5a_2 &2 &cgbc&  A_9^2D_6  &A_{12}^2   \cr
23 &24& &a_7^3 &12 &a&   A_8^3 & D_8^3  \cr
23 &25& &d_7a_4^4 &4 &bd&      &      \cr
23 &25& &a_9a_6a_4a_3 &1 &cfgb&      &      \cr
23 &25& &a_8d_5^2a_2^2 &2 &ebe&      &      \cr
23 &25& &d_6a_6a_6a_4 &1 &accf&      &      \cr
23 &25& &d_6a_6^2a_4 &2 &bce&      &      \cr
23 &25& &a_8a_7d_4a_3 &1 &fbbb&      &      \cr
23 &25& &a_7^2d_5a_3 &2 &baa&      &      \cr
23 &25& &a_6^2d_5^2 &2 &cb&      &      \cr
\cr
24 &20& &a_7^2d_5d_5 &4 &cda&   A_7^2D_5^2 &      \cr
24 &21& &a_8^2d_4a_4 &4 &baa&   A_8^3 &      \cr
24 &22& &a_9a_7d_4a_3a_1 &2 &fffad&  A_9^2D_6  &      \cr
24 &22& &a_7^2d_6a_3 &4 &cfa&  A_9^2D_6  &      \cr
24 &22& &d_6^2d_4^2a_3a_1^2 &4 &cbdb&  D_6^4 &      \cr
24 &24& &d_7a_5^2a_5a_1 &4 &cdea&  A_{11}D_7E_6 & A_{11}D_7E_6  \cr
24 &24& &a_9d_5^2a_1^2a_1 &4 &efec&  A_{11}D_7E_6 &A_{11}D_7E_6   \cr
24 &24& &a_7e_6d_4a_3^2 &4 &bgce&   A_{11}D_7E_6& A_{11}D_7E_6  \cr
24 &24& &e_6a_6^2a_4 &4 &eaa&  A_{11}D_7E_6 & A_{11}D_7E_6  \cr
24 &24& &e_6d_5a_5^2a_1^2 &4 &cbca& E_6^4 & A_{11}D_7E_6  \cr
24 &24& &d_6^2d_4^2a_1^4 &4 &cbb& D_6^4  &  D_8^3 \cr
24 &24& &a_7^2d_6a_1^2 &4 &dfd&  A_9^2D_6  &  D_8^3 \cr
24 &24& &a_7^2d_5d_4 &4 &fde&    A_7^2D_5^2&  A_{15}D_9 \cr
24 &25& &d_7a_5^2d_4a_1a_1 &2 &bdeeb&      &      \cr
24 &25& &a_9a_7a_3^2a_1 &2 &fgcd&      &      \cr
24 &25& &d_6d_5^2a_5a_1 &2 &bcba&      &      \cr
24 &25& &a_7d_6d_5a_3a_1a_1 &1 &cedeeb&      &      \cr
24 &25& &a_8a_7a_6a_1 &1 &aebd&      &      \cr
24 &25& &a_7d_6a_5d_4a_1 &1 &cfdfb&      &      \cr
24 &25& &a_7^2a_7a_1^2 &2 &edd&      &      \cr
24 &25& &a_8a_7d_4a_4 &1 &bfgb&      &      \cr
24 &25& &a_7a_7d_5a_4 &1 &ccgb&      &      \cr
24 &25& &a_7^2d_5a_4 &2 &cea&      &      \cr
\cr
25 &24& &a_{10}a_6a_5a_1 &2 &hgbb&  A_{11}D_7E_6 &  A_{12}^2 \cr
25 &24& &a_8a_7^2 &4 &eb&    A_8^3&  A_{15}D_9 \cr
25 &24&E &a_7^2d_5^2 &8 &aa&  2A_7^2D_5^2 &      \cr
25 &25& &d_7a_6a_6a_2a_2 &1 &adedc&      &      \cr
25 &25& &a_{10}a_6a_4a_2a_1 &1 &hgcea&      &      \cr
25 &25& &e_6a_6d_5a_4a_2 &1 &acaea&      &      \cr
25 &25& &e_6a_6^2d_4 &2 &bca&      &      \cr
25 &25& &a_7e_6a_5a_4a_1 &1 &acaeb&      &      \cr
25 &25& &a_{10}a_5a_4a_4 &1 &gbcb&      &      \cr
25 &25& &a_8d_6a_6a_2 &1 &dbfc&      &      \cr
25 &25& &a_9a_6d_5a_2a_1 &1 &bfbeb&      &      \cr
25 &25& &a_9a_6^2a_2 &2 &agc&      &      \cr
25 &25& &a_9d_5a_5a_4 &1 &bbbe&      &      \cr
\cr
26 &19& &a_7^2d_6d_5 &4 &dae&  A_7^2D_5^2  &      \cr
26 &21& &a_9d_6a_5d_4 &2 &cfea& A_9^2D_6   &      \cr
26 &23& &a_9d_6d_5a_1^2a_1 &2 &cffab& A_{11}D_7E_6  &      \cr
26 &23& &a_{10}a_6d_5a_1^2 &2 &bbga&  A_{11}D_7E_6 &      \cr
26 &23& &a_8e_6a_6a_2a_1^2 &2 &ahaba& A_{11}D_7E_6  &      \cr
26 &23& &d_7a_7d_5a_3a_1^2 &2 &edeea&  A_{11}D_7E_6 &      \cr
26 &23& &e_6d_6a_5^2a_1^2 &4 &dbea&  A_{11}D_7E_6 &      \cr
26 &23& &e_6d_5^3a_1^2 &12 &bce&  E_6^4&      \cr
26 &24& &a_9^2a_2^2 &8 &ga&  A_{12}^2 & A_{12}^2  \cr
26 &24& &d_7a_7d_5a_3 &2 &eddc&  A_{11}D_7E_6 &D_8^3   \cr
26 &24& &a_9d_6a_5a_3 &2 &dfbd&  A_9^2D_6  & A_{15}D_9  \cr
26 &24& &a_9a_8a_5 &2 &ebc&  A_9^2D_6  &  A_{15}D_9 \cr
26 &24& &a_7^2d_5^2 &8 &ce&  A_7^2D_5^2  &  D_{10}E_7^2\cr
26 &25& &d_7d_5^2a_3a_3 &2 &bdba&      &      \cr
26 &25& &a_{10}d_5a_5a_2a_1 &1 &bgbbb&      &      \cr
26 &25& &a_7d_6^2a_3 &2 &bcc&      &      \cr
26 &25& &a_9d_6d_4a_3a_1 &1 &cfbeb&      &      \cr
26 &25& &a_9a_7a_6 &1 &egb&      &      \cr
26 &25& &a_9a_7a_6 &1 &efb&      &      \cr
26 &25& &a_7d_6d_5a_5a_1 &1 &dffeb&      &      \cr
\cr
27 &24& &a_8^2a_7 &4 &ga&  A_8^3  &  A_{17}E_7 \cr
27 &25& &a_8d_7a_4a_4 &1 &dbde&      &      \cr
27 &25& &a_9^2a_3a_1^2 &2 &baa&      &      \cr
27 &25& &a_{10}a_6^2a_1 &2 &eca&      &      \cr
\cr
28 &20& &a_9^2d_5a_1^2 &4 &fac&  A_9^2D_6  &      \cr
28 &20& &d_6^3d_5a_1^2 &6 &bdb&  D_6^4 &      \cr
28 &22& &a_9e_6d_5a_3a_1 &2 &cfead& A_{11}D_7E_6  &      \cr
28 &22& &a_9d_7a_5a_3 &2 &dfca&   A_{11}D_7E_6&      \cr
28 &22& &a_{11}d_5a_5a_3a_1 &2 &fbeae&  A_{11}D_7E_6 &      \cr
28 &22& &e_6^2a_5^2a_3 &8 &bae&  E_6^4&      \cr
28 &23& &a_{10}a_9a_2a_1^2a_1 &2 &bgaaf& A_{12}^2  &      \cr
28 &23& &a_{11}a_7a_4a_1^2 &2 &gcaa& A_{12}^2  &      \cr
28 &24& &d_8d_4^4 &8 &cb&  D_8^3 &D_8^3   \cr
28 &24& &a_9d_7a_5a_1a_1 &2 &efede& A_{11}D_7E_6  & A_{15}D_9  \cr
28 &24& &a_{11}d_5d_4a_3 &2 &fbeb& A_{11}D_7E_6  & A_{15}D_9  \cr
28 &24& &a_9a_7d_6a_1 &2 &fefc&  A_9^2D_6  &  A_{17}E_7 \cr
28 &24& &a_9a_7d_6a_1 &2 &fbfc&  A_9^2D_6  & D_{10}E_7^2 \cr
28 &24& &d_6^2d_6d_4a_1^2 &2 &bbbb& D_6^4  &D_{10}E_7^2  \cr
28 &24&E &a_8^3 &12 &a& 2A_8^3  &      \cr
28 &25& &d_8a_5^2a_5a_1 &2 &ddbe&      &      \cr
28 &25& &a_{11}a_7a_2^2a_1^2 &2 &gbaf&      &      \cr
28 &25& &a_{11}d_5a_4a_3 &1 &fcbb&      &      \cr
28 &25& &e_6d_6d_5a_5a_1 &1 &cceab&      &      \cr
28 &25& &d_7d_6a_5a_5 &1 &cdcb&      &      \cr
28 &25& &a_8a_7e_6a_1a_1 &1 &aeedc&      &      \cr
28 &25& &a_9e_6a_4^2a_1 &2 &cgbc&      &      \cr
28 &25& &a_{10}a_8a_4a_1 &1 &bbae&      &      \cr
28 &25& &a_9^2a_4a_1^2 &2 &gaf&      &      \cr
28 &25& &a_9a_8d_5a_1 &1 &fbcd&      &      \cr
\cr
29 &24& &a_{11}a_8a_3 &2 &bca&  A_{12}^2 &  A_{15}D_9 \cr
29 &25& &a_{10}d_6a_6 &1 &fbe&      &      \cr
\cr
30 &21& &d_7a_7e_6d_4 &2 &cada&  A_{11}D_7E_6 &      \cr
30 &22& &a_{10}^2a_3 &4 &ba& A_{12}^2  &      \cr
30 &23& &d_7^2a_7a_1^2 &4 &dcd& D_8^3  &      \cr
30 &23& &d_8a_7^2a_1^2 &4 &ddd&  D_8^3 &      \cr
30 &24& &d_8a_7^2 &4 &dd&  D_8^3 & A_{15}D_9  \cr
30 &24& &a_{11}d_6a_5a_1 &2 &fbba& A_{11}D_7E_6  &A_{17}E_7   \cr
30 &24& &a_{10}e_6a_6 &2 &agb& A_{11}D_7E_6  & A_{17}E_7  \cr
30 &24& &d_7a_7e_6a_3 &2 &caed& A_{11}D_7E_6  & D_{10}E_7^2 \cr
30 &24& &a_9e_6d_6a_1 &2 &efea& A_{11}D_7E_6  & D_{10}E_7^2 \cr
30 &24& &e_6^2d_5^2 &8 &ca&  E_6^4& D_{10}E_7^2 \cr
30 &25& &d_8a_7d_5a_3 &1 &cbdd&      &      \cr
30 &25& &a_{10}a_9a_4 &1 &bca&      &      \cr
30 &25& &d_7a_7d_5^2 &2 &cae&      &      \cr
\cr
31 &24& &a_{12}a_7a_4 &2 &gbf& A_{12}^2  & A_{17}E_7  \cr
31 &24&E &d_6^4 &24 &d& 2D_6^4 &      \cr
31 &24&E &a_9^2d_6 &4 &aa& 2A_9^2D_6  &      \cr
31 &25& &a_{12}a_6a_5 &1 &gba&      &      \cr
31 &25& &a_8e_6^2a_2^2 &2 &aaa&      &      \cr
31 &25& &a_{10}e_6d_5a_2 &1 &ebba&      &      \cr
31 &25& &a_{11}a_8d_4 &1 &bea&      &      \cr
31 &25& &a_{11}a_7d_5 &1 &bba&      &      \cr
31 &25& &a_8a_8d_7 &1 &dcb&      &      \cr
31 &25& &a_8^2d_7 &2 &ca&      &      \cr
\cr
32 &18& &a_9^2d_7 &4 &da&  A_9^2D_6  &      \cr
32 &18& &d_7d_6^3 &6 &db& D_6^4  &      \cr
32 &20& &a_{11}e_6d_5a_3 &2 &ebaa& A_{11}D_7E_6  &      \cr
32 &21& &a_{12}a_8d_4 &2 &bba&  A_{12}^2 &      \cr
32 &22& &d_8d_6^2a_3a_1^2 &2 &bbdb&D_8^3   &      \cr
32 &24& &d_8d_6^2a_1^2a_1^2 &2 &bbba&D_8^3   &D_{10}E_7^2  \cr
32 &24& &d_6^4 &8 &b&  D_6^4 & D_{12}^2 \cr
32 &25& &a_9d_8a_5a_1a_1 &1 &dfecb&      &      \cr
32 &25& &a_{12}a_7d_4 &1 &bbf&      &      \cr
32 &25& &a_9d_7d_6a_1 &1 &ceeb&      &      \cr
\cr
34 &19& &a_{11}d_7d_6 &2 &cea&  A_{11}D_7E_6 &      \cr
34 &19& &e_6^3d_6 &12 &ae& E_6^4 &      \cr
34 &23& &a_{11}d_8a_3a_1^2 &2 &dbca& A_{15}D_9  &      \cr
34 &23& &a_{13}a_8a_1^2a_1 &2 &caac& A_{15}D_9  &      \cr
34 &23& &d_9a_7^2a_1^2 &4 &eea&  A_{15}D_9 &      \cr
34 &24& &a_{13}d_6a_3a_1 &2 &bbcb&A_{15}D_9   &A_{17}E_7   \cr
34 &24& &d_9a_7^2 &4 &ed& A_{15}D_9  & D_{10}E_7^2 \cr
34 &24& &a_{11}d_7d_5 &2 &eee& A_{11}D_7E_6  & D_{12}^2 \cr
34 &25& &e_7a_7d_5a_5 &1 &cbca&      &      \cr
34 &25& &e_7a_7^2a_3a_1 &2 &dcab&      &      \cr
34 &25& &d_9a_7d_5a_3 &1 &ddeb&      &      \cr
34 &25& &a_{13}a_7a_3a_1 &1 &cfcc&      &      \cr
34 &25& &d_7^2e_6a_3 &2 &aca&      &      \cr
34 &25& &d_8a_7e_6a_3 &1 &edda&      &      \cr
34 &25& &a_{11}d_6^2 &2 &cb&      &      \cr
\cr
35 &24& &a_{11}^2 &4 &a& A_{12}^2  & D_{12}^2 \cr
35 &25& &a_{12}d_7a_4 &1 &ebc&      &      \cr
\cr
36 &20& &d_8^2d_5d_4 &2 &bda& D_8^3  &      \cr
36 &22& &a_{13}d_7a_3a_1 &2 &ebab& A_{15}D_9  &   \cr
36 &24& &a_9e_7a_7 &2 &cfa&  D_{10}E_7^2& A_{17}E_7  \cr
36 &24& &e_7d_6d_6d_4a_1 &2 &bbbaa&D_{10}E_7^2  & D_{10}E_7^2 \cr
36 &24& &d_8^2d_4^2 &2 &bb&  D_8^3 &  D_{12}^2\cr
36 &25& &a_{12}a_{10}a_1 &1 &aab&      &      \cr
\cr
37 &24&E &e_6^4 &48 &e& 2E_6^4&      \cr
37 &24&E &a_{11}d_7e_6 &2 &aaa& 2A_{11}D_7E_6 &      \cr
37 &25& &a_{13}e_6a_4a_1 &1 &baba&      &      \cr
\cr
38 &17& &a_{11}d_8e_6 &2 &dae&  A_{11}D_7E_6 &      \cr
38 &21& &a_{11}d_9d_4 &2 &dea&  A_{15}D_9 &      \cr
38 &23& &d_9e_6^2a_1^2 &4 &add& D_{10}E_7^2 &      \cr
38 &23& &a_9e_7e_6a_1^2 &2 &aecd& D_{10}E_7^2 &      \cr
38 &23& &a_{14}e_6a_2a_1^2 &2 &bfaa&  A_{17}E_7 &      \cr
38 &23& &a_{11}e_7a_5a_1^2 &2 &cbba& A_{17}E_7  &      \cr
38 &24& &a_{11}d_9a_3 &2 &eeb&  A_{15}D_9 &  D_{12}^2\cr
38 &24& &a_{12}a_{11} &2 &af&  A_{12}^2 & A_{24}  \cr
38 &25& &d_9d_7a_7 &1 &cdb&      &      \cr
38 &25& &a_{11}d_8d_5 &1 &dbc&      &      \cr
\cr
40 &20& &a_{15}d_5d_5 &2 &bba&  A_{15}D_9 &      \cr
40 &22& &d_8e_7d_6a_3a_1 &1 &bbada&D_{10}E_7^2  &      \cr
40 &22& &d_{10}d_6^2a_3 &2 &bbd& D_{10}E_7^2 &      \cr
40 &22& &a_{15}d_6a_3 &2 &bca& A_{17}E_7  &      \cr
40 &24& &d_{10}d_6^2a_1^2 &2 &bba&D_{10}E_7^2  &D_{12}^2  \cr
40 &24&E &a_{12}^2 &4 &a&  2A_{12}^2&      \cr
40 &25& &d_{10}a_9a_5 &1 &dbb&      &      \cr
40 &25& &a_{15}d_5a_4 &1 &bca&      &      \cr
40 &25& &e_7e_6^2a_5 &2 &aaa&      &      \cr
40 &25& &a_9e_7d_7 &1 &aca&      &      \cr
40 &25& &a_{11}e_7d_5a_1 &1 &aeba&      &      \cr
40 &25& &a_{14}a_9 &1 &ac&      &      \cr
\cr
41 &24& &a_{15}a_8 &2 &ac& A_{15}D_9  &A_{24}   \cr
\cr
42 &21& &a_{13}e_7d_4 &2 &aba& A_{17}E_7  &      \cr
\cr
43 &24& &a_{16}a_7 &2 &fa&  A_{17}E_7 & A_{24}  \cr
43 &24&E &d_8^3 &6 &d&  2D_8^3&      \cr
\cr
44 &16& &d_9d_8^2 &2 &db& D_8^3  &      \cr
44 &20& &e_7^2d_6d_5 &2 &bad& D_{10}E_7^2 &      \cr
44 &24& &d_8^2d_8 &2 &ba&  D_8^3 & D_{16}E_8 \cr
44 &24& &d_8^3 &6 &a&  D_8^3 & E_8^3 \cr
\cr
46 &23& &d_{11}a_{11}a_1^2 &2 &ebd& D_{12}^2 &      \cr
46 &24& &a_{15}d_8 &2 &cb& A_{15}D_9  &D_{16}E_8  \cr
46 &25& &d_{11}a_7e_6 &1 &dab&      &      \cr
\cr
47 &25& &a_{16}d_7 &1 &ca&      &      \cr
\cr
48 &18& &d_{10}e_7d_7a_1 &1 &abda&  D_{10}E_7^2&      \cr
48 &18& &a_{17}d_7a_1 &2 &cac&  A_{17}E_7 &      \cr
48 &22& &d_{10}^2a_3a_1^2 &2 &bdb& D_{12}^2 &      \cr
48 &24& &a_{15}e_7a_1 &2 &bcc&  A_{17}E_7 & D_{16}E_8 \cr
48 &24& &d_{10}e_7d_6a_1 &1 &abaa& D_{10}E_7^2 & D_{16}E_8 \cr
48 &24& &d_8e_7^2a_1^2 &2 &aaa&  D_{10}E_7^2& D_{16}E_8 \cr
48 &25& &a_{13}d_{10}a_1 &1 &bcb&      &      \cr
\cr
49 &24&E &a_{15}d_9 &2 &aa& 2A_{15}D_9 &      \cr
\cr
50 &15& &a_{15}d_{10} &2 &ba& A_{15}D_9  &      \cr
\cr
52 &20& &d_{12}d_8d_5 &1 &bad&  D_{12}^2&      \cr
52 &23& &a_{19}a_4a_1^2 &2 &caa& A_{24}  &      \cr
52 &24& &d_{12}d_8d_4 &1 &bab&  D_{12}^2& D_{16}E_8 \cr
\cr
55 &24&E &d_{10}e_7^2 &2 &dd& 2D_{10}E_7^2&      \cr
55 &24&E &a_{17}e_7 &2 &aa& 2A_{17}E_7 &      \cr
55 &25& &a_{19}d_5 &1 &aa&      &      \cr
\cr
56 &14& &d_{11}e_7^2 &2 &da&  D_{10}E_7^2&      \cr
56 &21& &a_{20}d_4 &2 &aa&  A_{24} &      \cr
\cr
58 &25& &e_8a_{11}e_6 &1 &aac&      &      \cr
58 &25& &a_{11}d_{13} &1 &bb&      &      \cr
\cr
60 &24& &e_8d_8^2 &2 &aa& E_8^3 &D_{16}E_8  \cr
\cr
62 &23& &a_{15}e_8a_1^2 &2 &cbd&D_{16}E_8  &      \cr
\cr
64 &22& &e_8e_7^2a_3 &2 &aae& E_8^3 &      \cr
64 &22& &d_{14}e_7a_3a_1 &1 &abda& D_{16}E_8 &      \cr
\cr
67 &24&E &d_{12}^2 &2 &d& 2D_{12}^2&      \cr
\cr
68 &12& &d_{13}d_{12} &1 &db& D_{12}^2 &      \cr
68 &20& &d_{12}e_8d_5 &1 &aad& D_{16}E_8 &      \cr
68 &24& &d_{12}^2 &2 &b& D_{12}^2 & D_{24} \cr
\cr
71 &24& &a_{23} &2 &a& A_{24}  &  D_{24}\cr
\cr
76 &16&E&d_{16}d_9 &1 &bd& D_{16}E_8 &      \cr
76 &16&E&d_9e_8^2 &2 &ea& E_8^3 &      \cr
76 &24& &d_{16}d_8 &1 &ba&D_{16}E_8  &D_{24}  \cr
76 &24&E &a_{24} &2 &a&2A_{24}  &      \cr
\cr
91 &24&E &e_8^3 &6 &e& 2E_8^3&      \cr
91 &24&E &d_{16}e_8 &1 &dd&2D_{16}E_8 &      \cr
\cr
92 & 8&E&d_{17}e_8 &1 &da&D_{16}E_8  &      \cr
\cr
100 &20& &d_{20}d_5 &1 &ad&D_{24}  &      \cr
\cr
139 &24&E &d_{24} &1 &d& 2D_{24}&      \cr
\cr
140 & 0&E&d_{25} &1 &d& D_{24} &      \cr
}

\proclaim References.

\item{[A]}
Allcock, Daniel.
Recognizing Equivalence of Vectors in
    the Leech Lattice.
         Preprint 1996, available from 
\hbox{\tt  http://www.math.utah.edu/\~{}allcock}
\item{[B-V]}
Bacher, R.;  Venkov, B. B.
Reseaux entiers unimodulaires sans racines en dimension 27 et 28, 
Preprint 322, Inst. Fourier, Grenoble 1996. Available from
\hbox{\tt http://www-fourier.ujf-grenoble.fr/PREP/html/a332/a332.html}
\item{[B]}
Borcherds, Richard E.  The Leech lattice and other lattices, 1985 PhD
thesis, preprint math.NT/9911195 
\item{[B85]}
Borcherds, R. E. The Leech
lattice. Proc. Roy. Soc. London Ser. A 398 (1985), no. 1815, 365--376.
\item{[B90]}
 Borcherds, Richard E. Lattices like the Leech lattice. J. Algebra 130
 (1990), no. 1, 219--234.
\item{[C85]}
Conway, J. H. The automorphism group of the $26$-dimensional even
unimodular Lorentzian lattice. J. Algebra 80 (1983), no. 1, 159--163
\item{[C-P-S]}
Conway, J. H.; Parker, R. A.; Sloane, N. J. A. The covering radius of
the Leech lattice. Proc. Roy. Soc. London Ser. A 380 (1982), no. 1779,
261--290
\item{[C-S]}
Conway, J. H.; Sloane, N. J. A. Sphere packings, lattices and
groups. Third edition.  Grundlehren der Mathematischen
Wissenschaften 
290. Springer-Verlag, New York, 1999. ISBN: 0-387-98585-9
\item{[C-S98]}
Conway, J. H.;  Sloane, N. J. A. A Note on Optimal Unimodular Lattices 
J. Number Theory, 72 (1998), pp. 357--362. 
\item{[E]}
Eichler, M. ``Quadratische Formen und Orthogonal gruppen'', 
Springer--Verlag, 1952. 
\item{[E-G]}
Elkies, Noam D.; Gross, Benedict H. The exceptional cone and the Leech
lattice. Internat. Math. Res. Notices 1996, no. 14, 665--698.
\item{[K]}
Kneser, Martin. Klassenzahlen definiter quadratischer Formen.
Arch. Math. 8 (1957), 241--250.
\item{[N]}
Niemeier, H.-V. Definite quadratische Formen der Dimension 24 und
Diskriminante 1, J. Number theorey 5 (1973), 142--178.
\item{[S-V]}
Scharlau, Rudolf; Venkov, Boris B. The genus of the Barnes-Wall
lattice. Comment. Math. Helv. 69 (1994), no. 2, 322--333.
\item{[V75]}
Vinberg, \`E. B. 
Some arithmetical discrete groups in \hbox{Loba\v cevski\u\i}  spaces.
Discrete subgroups of Lie groups and applications to moduli
(Internat. Colloq., Bombay, 1973), pp.  323--348. Oxford Univ. Press,
Bombay, 1975.
\bye